\documentclass[12pt]{amsart}
\usepackage{amsfonts,amsmath,amssymb}
\usepackage{graphicx}
\usepackage[font=normalsize,textfont=normalfont]{caption}
\usepackage{amsmath,amssymb,amsfonts,amsthm,amstext,color,enumerate}
\usepackage{hyperref,cleveref,mathtools,mathdots}
\usepackage{mathtools}
\usepackage{url}
\usepackage{tikz}
\usetikzlibrary{arrows}

\newtheorem{maintheorem}{Theorem}

\newtheorem{lemma}[maintheorem]{Lemma}

\newtheorem{proposition}[maintheorem]{Proposition}
\newtheorem{corollary}[maintheorem]{Corollary}

\crefname{maintheorem}{Theorem}{Theorems}
\crefname{theorem}{Theorem}{Theorems}
\crefname{lemma}{Lemma}{Lemmas}
\crefname{proposition}{Proposition}{Propositions}
\crefname{corollary}{Corollary}{Corollaries}
\crefname{section}{Section}{Sections}
\crefname{figure}{Figure}{Figures}
\crefname{table}{Table}{Tables}

\renewcommand{\P}{\mathbb{P}}

\newcommand{\Z}{\mathbb{Z}}

\newcommand{\B}{\mathbb B}

\newcounter{mycount}
\newenvironment{ilist}{\begin{list}{\rm(\roman{mycount})}%
   {\usecounter{mycount}\labelwidth=8mm\itemsep 3pt}}{\end{list}}

\newcounter{acount}
\newenvironment{alist}{\begin{list}{\rm(\alph{acount})}%
   {\usecounter{acount}\labelwidth=8mm\itemsep 3pt}}{\end{list}}

\newcommand{\df}[1]{\textbf{#1}}

\begin{document}
\title{Games on Random Boards}
\author[Basu]{Riddhipratim Basu}
\address{Riddhipratim Basu, Department of Statistics, University of California,
Berkeley, USA}
\email{riddhipratim@stat.berkeley.edu}
\author[Holroyd]{Alexander E.~Holroyd}
\address{Alexander E.~Holroyd, Microsoft Research, Redmond, USA}
\email{holroyd@microsoft.com}
\author[Martin]{\\James B.~Martin}
\address{James B.~Martin, Statistics Department, University of Oxford, UK}
\email{martin@stats.ox.ac.uk}
\author[W\"astlund]{Johan W\"astlund}
\address{Johan W\"astlund, Department of Mathematical Sciences,
Chalmers University of Technology, Sweden}
\email{wastlund@chalmers.se}



\date{6 May 2015}
\begin{abstract}
We consider the following two-player game on a graph.  A
token is located at a vertex, and the players take turns
to move it along an edge to a vertex that has not been
visited before.  A player who cannot move loses.  We
analyse outcomes with optimal play on percolation
clusters of Euclidean lattices.

On $\Z^2$ with two different percolation parameters for
odd and even sites, we prove that the game has no draws
provided closed sites of one parity are sufficiently rare
compared with those of the other parity (thus favoring
one player).  We prove this also for certain
$d$-dimensional lattices with $d\geq 3$.  It is an open
question whether draws can occur when the two parameters
are equal.

On a finite ball of $\Z^2$, with only odd sites closed but with the external
boundary consisting of even sites, we identify up to logarithmic factors a
critical window for the trade-off between the size of the ball and the
percolation parameter. Outside this window, one or other
 player has a decisive advantage.

Our analysis of the game is intimately tied to the effect of boundary
conditions on maximum-cardinality matchings.
\end{abstract}

\maketitle

\section{Introduction}
Consider the following natural two-player game on an undirected graph. A
token is located at a vertex, and the players take turns to move. A move
consists of moving the token along an incident edge to a new vertex that has
never been visited before by the token. If a player has no possible move, she
loses (and the other player wins). We call this game \df{Trap} (since the
goal is to trap one's opponent).

We are concerned with optimal play. Thus, a \df{strategy} for the first or
second player is a map that assigns a legal next move (where one exists) to
each position. (A position comprises a location of the token and a set of
visited vertices.) Given a graph and an \df{initial vertex} (at which the
token starts), we say that the game is a \df{win} for the first or second
player respectively if that player has a \df{winning} strategy, i.e.\ a
strategy that results in a win, no matter what strategy the other player
uses. If the graph is finite, it is easy to check that the game is a win for
exactly one player.  In an infinite graph, it is possible that neither player
has a winning strategy, in which case we say that the game is a \df{draw}
(with the interpretation that the game continues forever with optimal play).

As we shall see, the outcome of Trap is intimately tied to
the properties of maximum-cardinality matchings, and draws
relate to sensitivity of such matchings to boundary
conditions.  In \cite{wastlund}, results about minimum
weight matchings in edge-weighted graphs were derived from
analysis of a related game called \emph{Exploration}.  We
shall consider another related game, which we call
\df{Vicious Trap},
 in which a player,
after making a move, is allowed to destroy (i.e., delete from the graph) any
subset of the vertices that he could have just moved to. These games are also
related to the game of Slither and its variants, as studied in \cite{And74,
Gardner,cornell}.

We are interested in playing {Trap} on a percolation
cluster. Let $G$ be an infinite connected graph, let
$p\in[0,1]$, and let each vertex of $G$ be declared
\df{closed} with probability $p$, and otherwise \df{open},
independently for different vertices.  Consider {Trap} on
the subgraph of $G$ induced by the set of open vertices.
(Equivalently, we play on $G$ but with moves to closed
vertices forbidden.  For convenience, the game is declared
a first player win if the initial vertex is closed.)  We
emphasize that the random subgraph is assumed known to both
players when deciding on their strategies.  If $p$ exceeds
one minus the critical probability for site percolation on
$G$, then almost surely the subgraph has no infinite
components, and therefore the game cannot be a draw.  For
many graphs (including the hypercubic lattice $\Z^d$), the
latter conclusion can be extended to a strictly longer
interval of $p$ using the method of essential enhancements
\cite{aizenman-grimmett}.  On $\Z^d$ itself (i.e.\ on the
open subgraph with $p=0$), Trap is easily seen to be a
draw. The regime of small positive $p$ seems to be the most
interesting.

A fascinating open question is whether there exists $p>0$ for which Trap on
the open subgraph of $\Z^d$ (started from the origin) is a draw with positive
probability.  As we discuss later, simulations lend some support to a
negative answer in dimension $d=2$.  On the other hand, draws \textit{do}
occur on certain random trees \cite{game-tree}, while variants of the model
on \textit{directed} lattices exhibit draws in dimensions $d\geq 3$ but not
in $d=2$ \cite{game-directed}.

Suppose now that the graph $G$ is bipartite, and call vertices in its two
classes \df{odd} and \df{even}.  A natural extension of the above model is to
declare odd and even vertices closed with respective probabilities $p$ and
$q$ (with different vertices still receiving independent assignments). Given
an initial vertex (which may be odd or even), one player always moves from
even vertices to odd vertices. Call this player \df{Odin}, and the other
player \df{Eve}. We summarize these conventions in \cref{odd-even}. The
probability that Odin wins is non-increasing in $p$ and non-decreasing in
$q$, and vice-versa for Eve. (Indeed, introducing more closed odd vertices
preserves all winning strategies for Eve but cannot create winning strategies
for Odin.)  However, there is apparently no obvious monotonicity argument for
the probability of a draw.
\begin{table}
\renewcommand{\arraystretch}{1.2}
\makebox[\textwidth][c]{
 \begin{tabular}{|c|c|c|}
   \hline
  \textbf{\textcolor{red}{Odin}} moves to & \textbf{odd} vertices,
   & which are closed with probability $p$.   \\
  \hline
  \textbf{\textcolor{blue}{Eve}} moves to & \textbf{even} vertices,
   & which are closed with probability $q$. \\
   \hline
 \end{tabular}
}
 \caption{Conventions for bipartite graphs.}\label{odd-even}
\end{table}

\enlargethispage*{2\baselineskip}
 If $p>q$ then we should
expect Eve to have an advantage. Our first main result
states that this advantage is decisive in the extreme case
$q=0$.  We prove this in all dimensions $d$, but on a
slightly non-standard lattice when $d\geq 3$. The
\df{body-centered hypercubic lattice} is the graph $\B^d$
whose vertices are all elements of $\Z^d$ having
coordinates that are all even or all odd (called even and
odd vertices respectively), and with an edge between
vertices $u$ and $v$ whenever $\|u-v\|_\infty=1$. Note that
$\B^2$ is isomorphic to the usual square lattice $\Z^2$.

\begin{maintheorem}
\label{bcc}
Let $d\geq 2$ and consider the body-centered lattice
$\B^d$, with odd and even vertices closed with respective
probabilities $p$ and $q$. Consider a game of Trap between
 Eve and Odin on the open subgraph of $\B^d$.
\begin{ilist}
\item Let $p>0$ and $q=0$.  Almost surely, for every
    initial vertex, Eve wins.
\item For all $p>0$ there exists $q(p,d)>0$ such that if
    $q<q(p,d)$ then almost surely, for every initial vertex,
     the game is not a draw.
\end{ilist}
\end{maintheorem}

 The key step in our proof of \cref{bcc} will be to
show the existence of finite regions with the property that if Odin enters
one, he cannot escape. This is proved using a result of \cite{schonmann} on
(modified) bootstrap percolation. To ensure that with high probability a
region has the required properties for this argument, it must be very large,
of order $\exp^{d-1} (\lambda/p)$ for small $p$ (where $\lambda=\pi^2/6$, and
the exponential function is iterated $d-1$ times). This follows from results
of \cite{cerf-cirillo,cerf-manzo,holroyd-2d,holroyd}. The resulting lower
bound on $q(p,d)$ is therefore very small: $O(1/\exp^{d-1} (\lambda/p))$ as
$p\to 0$.

It is likely that our methods could be adapted to prove
that the conclusions of \cref{bcc} hold also for the
standard hypercubic lattice $\Z^d$ in all dimensions.
Checking this entails adapting standard results of
\cite{schonmann} to a variant of bootstrap percolation in a
different combinatorial setting.  Rather than pursuing
this, we focus next on the more interesting question of
obtaining tighter bounds involving finite regions in
dimension $2$. This will yield improved bounds on $q(p,2)$.
Moreover, understanding the game on finite regions is an
important step towards the main open question about the
case $p=q$ as discussed earlier.

\begin{figure}
  \centering
  \begin{tabular}{cp{.05\textwidth}c}
     \includegraphics[width=.375\linewidth]{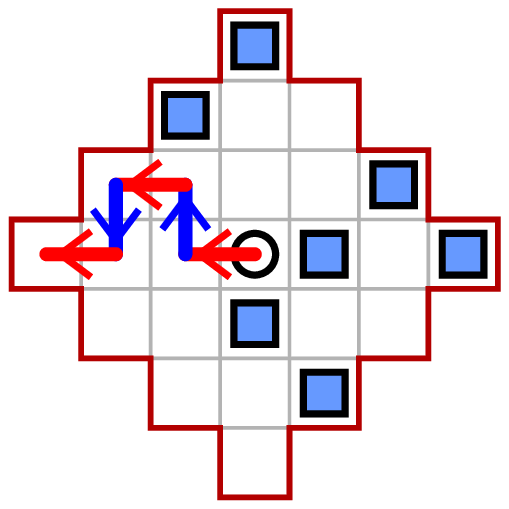}
     &&
     \includegraphics[width=.375\linewidth]{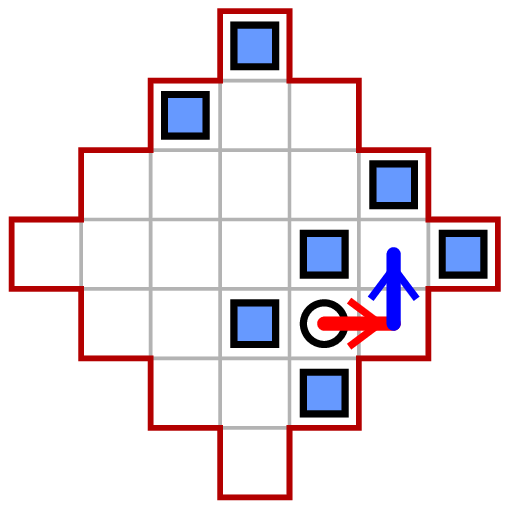}
     \\[-1mm]
     Odin starts and wins.
     &&
     Odin starts and loses. \\[3mm]
     \includegraphics[width=.375\linewidth]{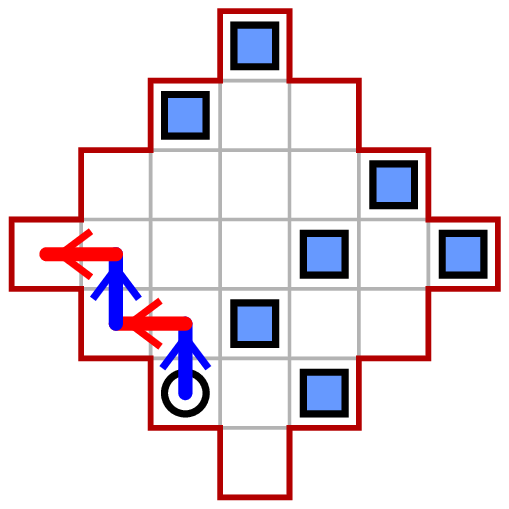}
     &&
     \includegraphics[width=.375\linewidth]{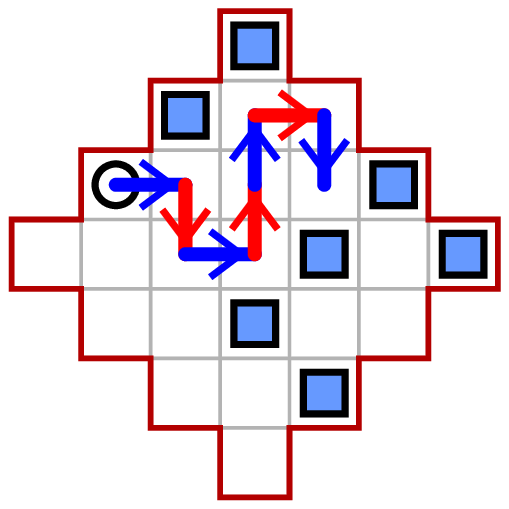}
     \\[-1mm]
     Eve starts and loses.
     &&
     Eve starts and wins. \\
  \end{tabular}
\caption{Examples of optimal play on the diamond $D_2$, from four different
initial vertices.  Squares represent vertices, with closed (odd) vertices
shown as filled squares, and the initial vertex marked with a circle. Odin's
moves are red; Eve's moves are blue.  (The colors of the boundary and the
closed vertices reflect the fact that they are helpful to Odin and Eve
respectively.)}\label{f:gameplay}
\end{figure}
%
%
%
We restrict attention to the square lattice $\Z^2$ (i.e.\
the graph with vertex set $\Z^2$, and an edge between $u$
and $v$ whenever $\|u-v\|_1=1$). A vertex is called
\df{odd} or \df{even} according to whether the sum of its
coordinates is odd or even.  We take $p>0$ and $q=0$ (so
that only odd vertices can be closed), but we offset Eve's
advantage by restricting to a finite region in a way that
favors Odin. Specifically, let $n>0$ be an integer, and let
$D_n$ be the subgraph of $\Z^2$ induced by the
 region $\{u\in\Z^2:\|u\|_1 <2n\}$. We call
$D_n$ a \df{diamond}. Note that all internal boundary
vertices of $D_n$ (i.e.\ those $u$ with $\|u\|_1=2n-1$) are
odd.  We consider the game on the open subgraph of $D_n$.
Thus, Odin is forbidden from moving to closed vertices, but
Eve is forbidden from moving out of $D_n$. Intuitively, Eve
tries to trap Odin using the closed odd vertices, while
Odin tries to trap Eve against the boundary. (Equivalently,
we can consider the game on the open vertices of $\Z^2$,
but declaring a win for Odin if the token ever leaves
$D_n$, or alternatively we can declare all vertices outside
$D_n$ to be closed). The progress of this game in a few
cases is illustrated in Figure \ref{f:gameplay}.

It is easy to see that the probability that Eve wins (starting from the
origin, say) is non-decreasing in $p$ and $n$. We address how these two
effects compare with each other as $(n,p)\to (\infty,0)$. We show that up to
logarithmic factors, the critical regime is at $n = \mbox{constant}/p$;
outside this window, one player has a decisive advantage, for essentially
every initial vertex. In the following, ``with high probability" means with
probability tending to $1$.

\begin{maintheorem}\label{lower}
Let odd and even vertices of the diamond $D_n$ be closed
with respective probabilities $p>0$ and $q=0$, and consider
a game of Trap. For every constant $c>0$, if $n< c /(p \log
p^{-1})$, then, with high probability as $p\to 0$, from
every initial vertex in $D_n$, Odin wins.
\end{maintheorem}

The next result gives complementary conditions under which
Eve wins, but now there are exceptional initial vertices
near the boundary.  Define $K_0=\{(x,y)\in \Z^2: |y|<x\}$,
and for $k=1,2,3$ let $K_k$ be obtained by rotating $K_{0}$
counter-clockwise by $\pi k/2$ about the origin.  We say
that a vertex $u\in D_n$ is \df{protected} if each of the
cones $u+K_0,\ldots,u+K_3$ contains some closed vertex of
$D_n$. If the initial vertex is even and unprotected, then
a simple winning strategy for Odin is to choose such a cone
containing no closed vertex, and always move in the
direction of that cone (e.g.\ rightwards in the case of
$u+K_0$).

\begin{samepage}
\begin{maintheorem}\label{upper}
Let odd and even vertices of the diamond $D_n$ be closed
with respective probabilities $p>0$ and $q=0$, and consider
a game of Trap. There exists a constant $C>0$ such that, if
$n> (C \log p^{-1} )/p$, then, with high probability as
$p\to 0$, from every odd vertex and every protected even
vertex of $D_n$, Eve wins.
\end{maintheorem}
\end{samepage}

It is straightforward to check that for a fixed constant $C'$, with high
probability, every even vertex in the set
\begin{equation}\label{hyperbola}
S:=\left\{(x,y)\in D_n: \Bigl(2n- \Bigl|\frac{x+y}{2}\Bigr|\Bigr)
\Bigl(2n-\Bigl|\frac{x-y}{2}\Bigr|\Bigr) > \frac{C'\log
p^{-1}}{p}\right\}
\end{equation}
 is protected. Since $|D_n\setminus S|= O(p^{-1}\log
^2 p^{-1})$ as $p\rightarrow 0$ (uniformly in $n$), in the
situation of \cref{upper}, the conclusion of the theorem
applies to all but a fraction $O(p)$ of the even vertices
of $D_n$ with high probability as $p\rightarrow 0$. (We
justify these remarks in \S~\ref{s:game}.)

One consequence of \cref{upper} is that for the $2$-dimensional lattice
$\Z^2$ (i.e.\ $\B^2$), \cref{bcc} (ii) holds with $q(p,2)=c' p^2\log ^{-2}
p^{-1}$ for some absolute constant $c'$, much better than the bound
$\exp(-C'' /p)$ that results from bootstrap percolation arguments.  See
\S~\ref{s:game} for details. Using \cref{lower,upper} we can also provide
upper and lower bounds for the number of steps required for the game to
terminate on $\Z^2$; see \S~\ref{s:game}.

\subsection*{Matchings and independent sets}

Our proofs of \cref{lower,upper} rely on the following characterization of
winning positions for {Trap} on a finite graph, which we prove in
\S~\ref{s:finite}. A closely related result appears in \cite{cornell}, for a
game that is a variant of Trap.

\begin{proposition}
\label{t:matchgame} Let $G=(V,E)$ be a finite, connected,
simple graph. Trap on $G$ starting from $v\in V$ is a win
for the first player if and only if $v$ is contained in all
maximum-cardinality matchings of $G$.
\end{proposition}

Now we consider Vicious Trap. Recall that this game is same as Trap except
that in addition to moving the token along an edge to a previously unvisited
vertex, a player may delete any subset of the vertices that he could have
moved to (i.e.\ of the neighbors of the previous vertex other than the
current vertex). Moves to deleted vertices are forbidden, and a player who
cannot move loses. For this game we provide an analogous characterization of
winning positions, now involving maximum-cardinality independent sets.

\begin{proposition}
\label{t:indgame} Let $G=(V,E)$ be a finite, connected, simple graph. Vicious
Trap on $G$ starting from $v\in V$ is a loss for the first player if and only
if $v$ is contained in every maximum-cardinality independent set of $G$.
\end{proposition}

{Trap} and {Vicious Trap} are equivalent on a bipartite graph $G$.  This is
because it is never advantageous for a player to delete vertices, since those
vertices are inaccessible to the other player anyway. By
\cref{t:matchgame,t:indgame}, this gives an interesting proof that on a
finite simple bipartite graph, a vertex is contained in all maximum matchings
if and only if it is absent from some maximum independent set. This fact can
be deduced more directly from {{K\"{o}nig's Theorem}} which states that on a
bipartite graph, the number of edges in a maximum matching equals the number
of vertices in a minimum vertex cover \cite[Theorem~5.3]{Konig}.

Since the graphs we consider ($d$-dimensional lattices and subgraphs thereof)
are bipartite, the conclusions of \cref{bcc,lower,upper} hold for {Vicious
Trap} as well.

\begin{figure}
  \centering
  \begin{tabular}{cc}
     \includegraphics[width=.49\linewidth]{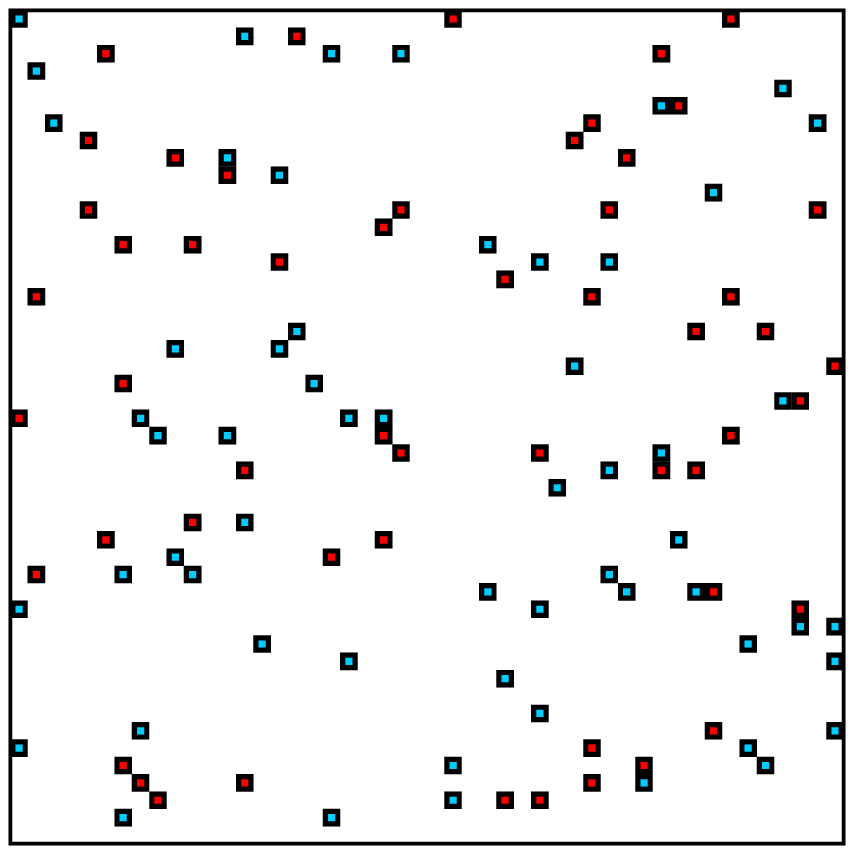}
     &
     \includegraphics[width=.49\linewidth]{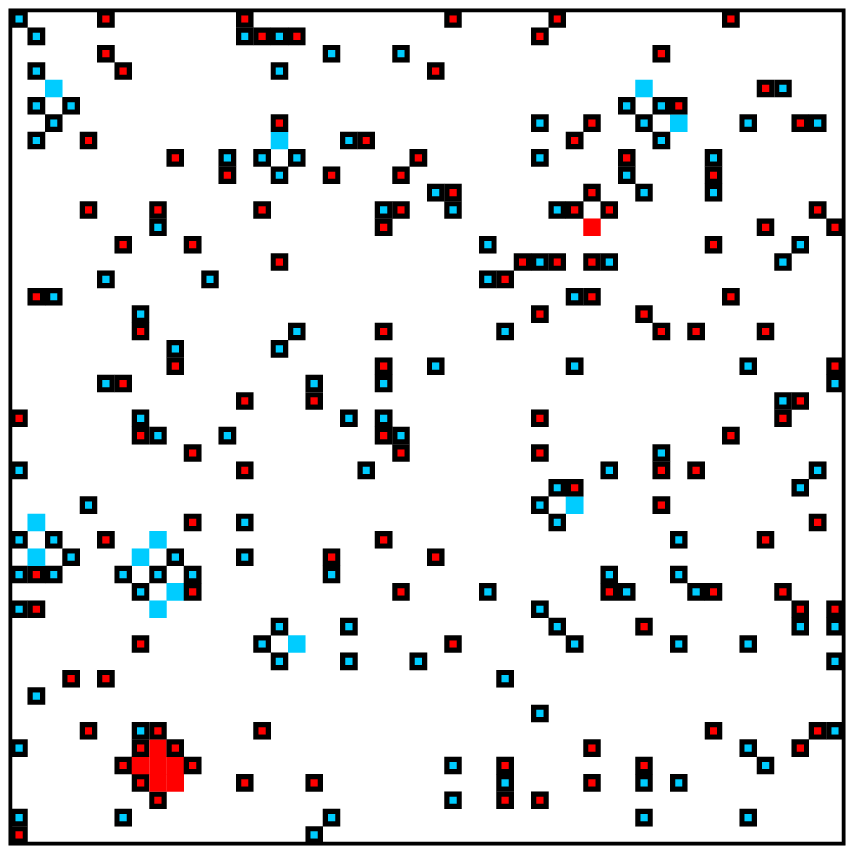}
     \\[-1mm]
     $p=0.05$
     &
     $p=0.1$ \\[3mm]
     \includegraphics[width=.49\linewidth]{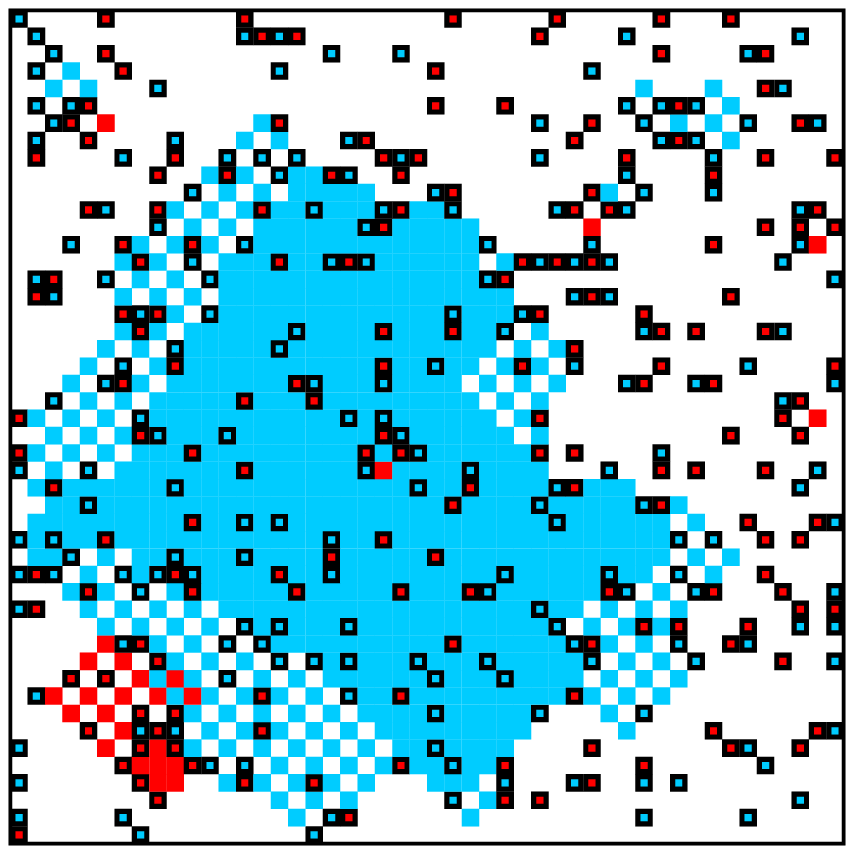}
     &
     \includegraphics[width=.49\linewidth]{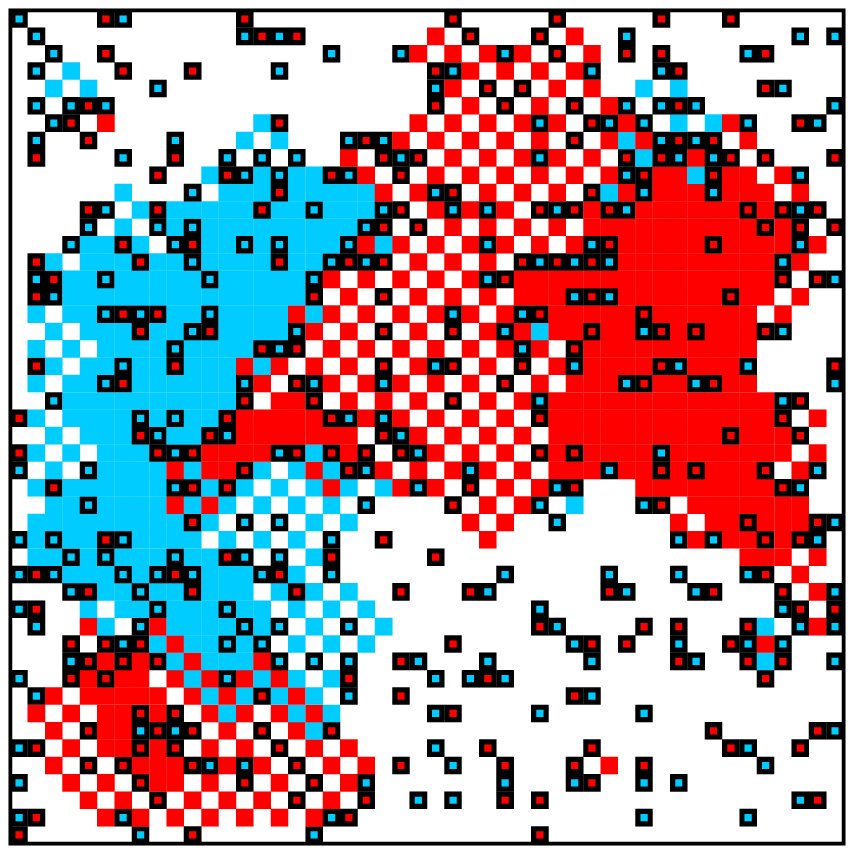}
     \\[-1mm]
     $p=0.15$
     &
     $p=0.2$ \\
  \end{tabular}
%
\caption{Outcomes of Trap on a square of size $n=50$, with
the game declared a draw if the token leaves the square.
Closed vertices occur with probability $p$, and are
outlined in black, with their interiors respectively blue
or red according to whether they are odd or even (to
reflect the fact that they are favorable to Eve or Odin
respectively). Other vertices are blue if Eve wins, red if
Odin wins, or white if the game is drawn, from that initial
vertex.} \label{f:simulation50}
\end{figure}
\begin{figure}
  \centering
  \begin{tabular}{cc}
     \includegraphics[width=.49\linewidth]{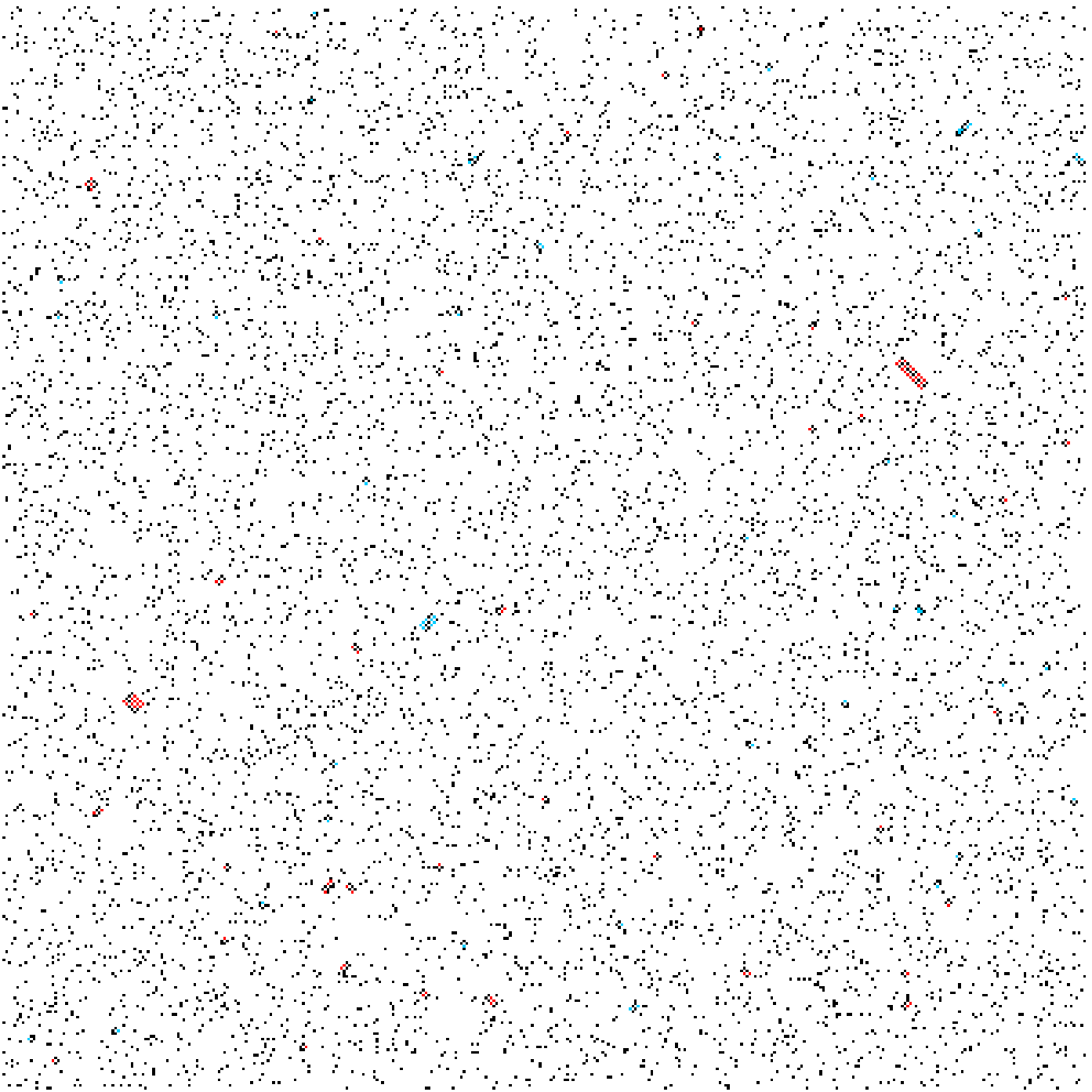}
     &
     \includegraphics[width=.49\linewidth]{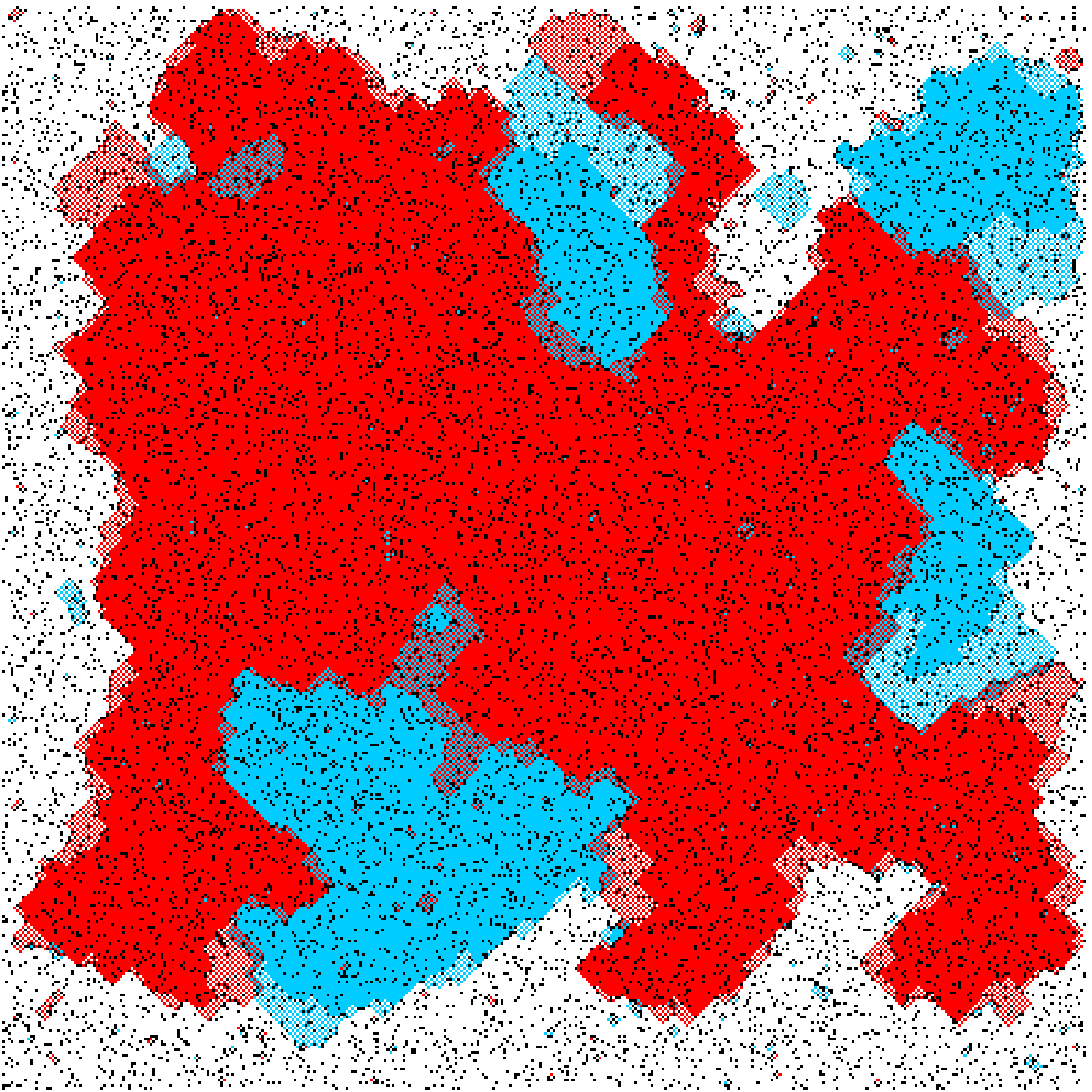}
     \\
     $p=0.05$
     &
     $p=0.1$ \\[3mm]
     \includegraphics[width=.49\linewidth]{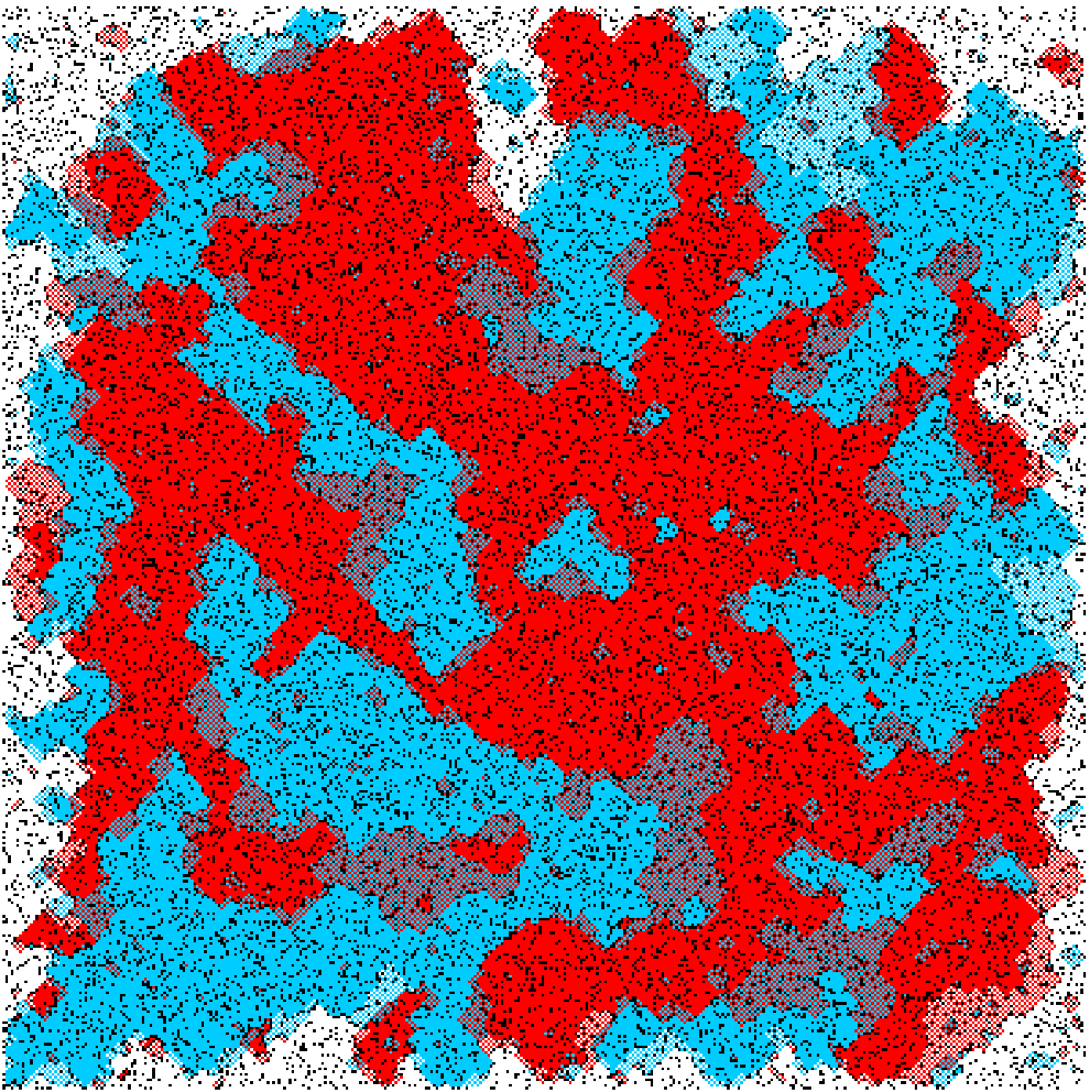}
     &
     \includegraphics[width=.49\linewidth]{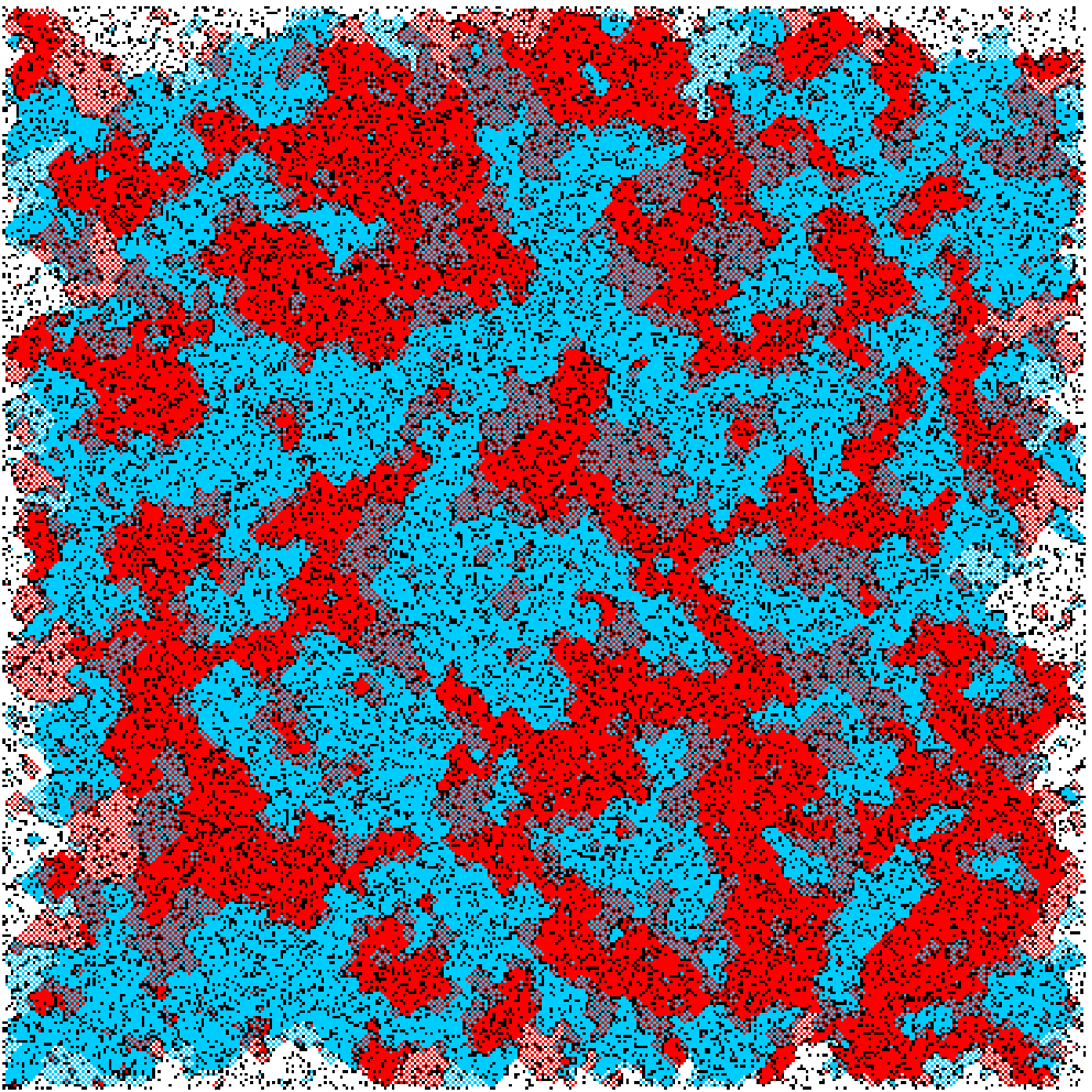}
     \\
     $p=0.15$
     &
     $p=0.2$ \\
  \end{tabular}
%
\caption{Outcomes of Trap by initial vertex on a square of
size $n=400$, with the boundary declared a draw. Closed
vertices are black, wins for Eve are blue, wins for Odin
are red, and draws are white.} \label{f:simulation400}
\end{figure}
\subsection*{Simulations and Conjectures}
\cref{t:matchgame} gives rise to a practical algorithm for
determining the outcome of Trap on a finite bipartite
graph.  We can find a maximum size matching $M$ using the
Hopcroft-Karp algorithm \cite{HK}.  Then we can search for
all matched vertices $v$ from which there is no alternating
path leading to an unmatched vertex.  A vertex $v$ has this
property precisely if it is contained in every maximum
matching, i.e.\ if it is a winning initial vertex for the
first player.

To gain insight about Trap in the most interesting setting of $\Z^2$ with
$p=q$, we may proceed as follows.  Consider the square $[1,n]^2\cap\Z^2$,
with odd and even vertices closed with equal probability $p=q$, and declare
the game a draw if the token ever leaves the square.  We may determine the
outcome of this game by applying the method described above to the square
with two different boundary conditions, and comparing the results.  In one
case we modify the graph just outside the boundary of the square so that all
internal boundary vertices of the resulting graph are even; this means that
Eve wins if the token leaves the square.  In the other case we similarly
arrange that Odin wins if the token leaves the square.  An initial vertex
should then be considered a draw if its outcomes differ between these two
boundary conditions. Thus we can identify the outcome from every initial
vertex.

\cref{f:simulation50,f:simulation400} show the results of the above
experiment on squares of sizes $50$ and $400$, with four different values of
$p$.  The results are suggestive of the following picture.  For each $p$,
there are domains within which one or other player can force a win,
presumably owing to a local preponderance of closed vertices of the
appropriate parity. These domains tend to abut each other, so that there are
no regions of draws between them.  (An additional complication is the
appearance of ``checkerboard'' regions, from which the first player wins,
near some interfaces between opposing domains).  The typical size of a domain
apparently diverges as $p$ becomes small, and, if the square is not large
enough to contain a whole domain, then instead draws are prevalent.

However, the simulations seem consistent with the
hypothesis that the typical domain size is finite for each
$p>0$, and only approaches $\infty$ as $p\to 0$.  This
would suggest that there are no draws on $\Z^2$ for any
$p>0$. (Other interpretations of the data are possible, and
our confidence in this conclusion is not especially high).

On the other hand, on $\Z^d$ with $d\geq 3$, by analogy
with the directed variants of Trap considered in
\cite{game-directed}, one may speculate that draws {\em do}
occur when $p$ and $q$ are equal and sufficiently small.

\subsection*{Organization of the paper}
The rest of this paper is organized as follows. In \S~\ref{s:bcc} we prove
\cref{bcc} using a bootstrap percolation argument. Next, in
\S~\ref{s:finite}, we prove \cref{t:matchgame,t:indgame}.  In \S~\ref{s:ub}
and \S~\ref{s:lb} respectively we use these results to prove Theorems
\ref{upper} and \ref{lower} on the diamond.  In \S~\ref{s:game}, we give some
consequences for regions of other shapes, and for the length of the game on
$\Z^2$.

\section{Bootstrap Percolation Bound}
\label{s:bcc}

Our proof of \cref{bcc} exploits a connection with a variant of bootstrap
percolation.  The basic idea is as follows. Let $u$ be an open even vertex of
the body-centered lattice $\B^d$, and suppose that of the $2^d$ adjacent odd
vertices, exactly one, $v$, is open. Starting from $v$, Eve can win
immediately by moving to $u$.  Thus, $v$ is effectively forbidden to Odin,
and so we can now iterate the argument with $v$ added to the set of closed
odd vertices.

Here is the relevant bootstrap percolation model on $\Z^d$,
which we call the \df{Fr\"obose model} because it is a
natural extension to $d$ dimensions of a model introduced
in \cite{Frobose}.  We start with a given subset $X_0$ of
$\Z^d$, whose elements are said to be \df{occupied} at time
$0$.  We define the set of occupied vertices at time $t$,
denoted $X_t$, for $t>0$, inductively as follows.  Any
vertex occupied at time $t-1$ remains occupied at time $t$.
In addition, if all but one of the elements of any
hypercube of the form $u+\{0,1\}^d$ are occupied at time
$t-1$, then the one remaining vertex of the hypercube
becomes occupied at time $t$.  Let $\langle
X_0\rangle=X_\infty:=\bigcup_{t=0}^\infty X_t$ be the set
of eventually occupied vertices.

Given an initially occupied set $X_0\subseteq \Z^d$, a set
$W\subseteq \Z^d$ is said to be \df{internally spanned} if
$W\subseteq \langle X_0\cap W\rangle$, i.e.\ if $W$ becomes
fully occupied when we start from only the initially
occupied vertices in $W$.

Let $B(n):=[1,n]^d\cap\Z^d$. The following is a standard result of bootstrap
percolation, adapted to the Fr\"obose model.
\begin{proposition}
\label{p:frobose} Fix $p$ and let each vertex of $\Z^d$ be
initially occupied independently with probability $p$.  For
all $d\geq 1$ and $p>0$ we have
$$\P\bigl(B(n)~\text{\rm is internally spanned}\bigr)\rightarrow 1
    \text{ as } n\rightarrow \infty.$$
\end{proposition}

Proposition \ref{p:frobose} follows from arguments of
\cite{schonmann}. Another version of the argument, giving
much tighter bounds, appears in \cite{holroyd}.  The
relevant results in \cite{schonmann,holroyd} state that the
conclusion of \cref{p:frobose} holds for another model
called {\em modified} bootstrap percolation (in which a
vertex becomes occupied if it has at least one occupied
neighbor in each dimension). The Fr\"obose model is
`weaker' in the sense that it is harder for a vertex to
become occupied, so the conclusion itself does not carry
over directly. However, the proofs in
\cite{schonmann,holroyd} proceed by defining particular
events $\mathcal{E}_n$ such that $\P(\mathcal{E}_n)\to 1$
as $n\to\infty$, and such that $B(n)$ is internally spanned
(with respect to the modified model) on the event
${\mathcal E}_n$.  It turns out (and it is straightforward
to verify) that $B(n)$ is also internally spanned with
respect to the Fr\"obose model on the same event
$\mathcal{E}_n$, so \cref{p:frobose} immediately follows.

To connect the Fr\"obose bootstrap model with Trap on
$\B^d$, let $\B_{\rm o}^d$ (respectively, $\B_{\rm e}^d$)
be the graph comprising all odd (even) vertices of $\B^d$,
with an edge between any pair of vertices that are at
$\ell_{1}$ distance $2$. Obviously $\B_{\rm o}^d$ is
isomorphic to the standard hypercubic lattice $\Z^d$.  Let
vertices of $\B^d$ be open or closed, as usual.  Then we
may run the Fr\"obose bootstrap percolation model on
$\B_{\rm o}^d$ (by which we mean that we consider the image
under the obvious isomorphism of the model on $\Z^d$).  We
declare the closed odd vertices initially occupied.

Recall that $B(n):=[1,n]^d\cap\Z^d$.  For an odd vertex
$u\in \B_{\rm o}^d$, define the set of odd vertices
$\widetilde{B}_{\rm o}(u,n):=u+2B(n) \subset \B_{\rm o}^d$.
Note that the induced subgraph of $\widetilde{B}_{\rm
o}(u,n)$ in the graph $\B_{\rm o}^d$ is isomorphic to the
induced subgraph of $B(n)$ in $\Z^d$. Also let
$\widetilde{B}_{\rm e}(u,n)\subset \B_{\rm e}^d$ be the set
of all even vertices all of whose neighbors lie in
$\widetilde{B}_{\rm o}(u,n)$, and let
$\widetilde{B}(u,n)=\widetilde{B}_{\rm o}(u,n)\cup
\widetilde{B}_{\rm e}(u,n)$.

We call the box $\widetilde{B}(u,n)$ \df{good} if (i) all of its even
vertices are open, and (ii) its odd subgraph $\widetilde{B}_{\rm o}(u,n)$ is
internally spanned with respect to the Fr\"obose bootstrap model on $\B_{\rm
o}^d$ started with the closed vertices occupied.

\begin{proposition}
\label{p:bccbox} Suppose that the box $\widetilde{B}(u,n)$ is good.  Then
from every odd $v\in \widetilde{B}(u,n)$, Eve has a winning strategy for Trap
that guarantees that the token never leaves $\widetilde{B}(u,n)$.
\end{proposition}

\begin{proof}
Let $X_0$ be the set of closed odd vertices in
$\widetilde{B}_{\rm o}(u,n)$, and let $X_t\subseteq
\widetilde{B}_{\rm o}(u,n)$ be the set of vertices that are
occupied at time $t$ starting from $X_0$ occupied.  For
$v\in \widetilde{B}_{\rm o}(u,n)$, let $T(v)$ be the time
at which $v$ becomes occupied, i.e.\ let $T(v)=t$ if $v\in
X_t\setminus X_{t-1}$.

For the purpose of this proof it is convenient to allow
Odin to move to a closed vertex, but declare an immediate
win for Eve if he does so.  (This clearly does not change
the outcome of the game).  We claim that starting from any
$v\in \widetilde{B}_{\rm o}(u,n)$, Eve has a winning
strategy that guarantees that $T$ is strictly decreasing
along the sequence of odd vertices that are visited.  Here
is Eve's strategy, which we define inductively. Suppose
(perhaps after some steps of the game played according to
such a strategy) that it is Eve's turn. Then the token is
at an odd vertex $v$; suppose $T(v)=t\geq 1$. By definition
of the Fr\"obose model, there exists an even neighbor $w$
of $v$, all of whose neighbors other than $v$ lie in
$X_{t-1}$. Since $v$ is the first vertex in $X_t$ that has
been visited, no other neighbor of $w$ has been visited
before, therefore $w$ has not been visited. Therefore, Eve
moves to $w$, and Odin must then move to an element of
$X_{t-1}$, as required.

In particular, we deduce that starting from a vertex in $X_t$, Eve has a
winning strategy that guarantees that the token remains in
$\widetilde{B}(u,n)$ and that Eve makes at most $t$ moves.
\end{proof}

Standard results also show that $\Z^d$ itself is internally spanned almost
surely in the Fr\"obose model for any positive density $p$ of initially
occupied sites.  Hence, a minor variant of the above argument already shows
that with $p>0$ and $q=0$, Eve wins on $\B^d$ if she has the first move. To
deal with small positive $q$ and the possibility that Odin starts, we need to
be a little more careful.

\begin{proof}[Proof of \cref{bcc}]
Let $0:=(0,\ldots, 0)\in\B^d_{\rm e}$ and
$\iota:=(1,\dots,1)\in\B^d_{\rm o}$.  By translation
invariance, it suffices to prove the claims for Trap
started at $0$ or $\iota$.  For an integer $n\geq 1$ and
$x\in \Z^d$ we introduce the renormalized box $\widehat B
(x,n):=\widetilde B(\iota+2n x,n)$.  Note that these boxes
are disjoint for different $x$, but adjacent boxes almost
abut each other: if $\|x-y\|_1=1$ then there is a layer of
even vertices between $\widehat B (x,n)$ and $\widehat B
(y,n)$, but all neighbors of those vertices lie in one of
the two boxes.

We first prove part (ii).  Let $p_{\rm c}$ be the critical
probability of site percolation on the $d$-dimensional
\df{star-lattice}, i.e.\ the graph with vertex set $\Z^d$
and an edge between $u$ and $v$ whenever
$\|u-v\|_{\infty}=1$. Fix $p>0$.  By \cref{p:bccbox}, there
exists $n\geq 1$ such that with $q=0$, we have
$\P(\widehat{B}(0,n)\text{ is good})>1-p_{\rm c}$.  Since
the box has finitely many even vertices, the same
conclusion holds for $q$ sufficiently small; fix such a
$q$.

\begin{figure*}
\begin{center}
\includegraphics[width=0.45\textwidth]{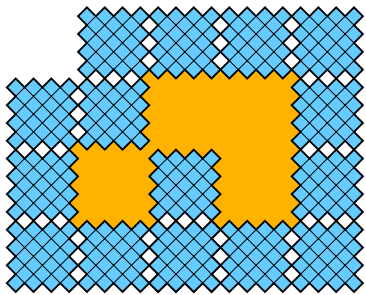}
\caption{A finite region $S$ surrounded by good boxes.  (Vertices of $\B^2$
are shown as squares at 45 degrees to the axes). Odin can only leave $S$ by
entering a good box, whereupon Eve can win. } \label{f:contour}
\end{center}
\end{figure*}
Let $Z:=\{x\in\Z^d: \widehat{B}(x,n) \text{ is not good}\}\cup\{0\}$.  Note
that we include the origin regardless of whether or not $\widehat B (0,n)$ is
good. The set $Z$ induces a subgraph of the star-lattice; let $K$ be the
vertex set of the component containing $0$. Then $K$ is finite almost surely.
Let $S$ be union of the renormalized boxes $\widehat{B}(x,n)$ for $x\in K$,
together with the set of all even vertices of $\B^d$ adjacent to them.  Then
$S$ is a.s.\ finite and contains $\widehat B (0,n)$ (which contains $0$ and
$\iota$). Furthermore, every infinite path in $\B^d$ starting from a vertex
in $S$ intersects some good box that is not $\widehat B (0,n)$.  See
\cref{f:contour}.

Suppose that the token starts from $0$ or $\iota$ and at
some point leaves $S$.  Then it must do so by entering a
good box, and it must enter it at an odd vertex, via a move
of Odin. \cref{p:bccbox} implies that Eve can then win
within that box.  It follows that the outcome of Trap
starting from $0$ or $\iota$ is identical to the outcome
restricted to $S$ (i.e.\ with moves out of $S$ forbidden,
and the same vertices closed as before). Since $S$ is
finite, the game cannot be a draw, thus establishing (ii).

Now we prove (i).  For all $p>0$, the conclusion of (ii)
applies when $q=0$, so the game is not a draw.  However,
Eve has a simple strategy that guarantees she cannot lose:
she always moves in direction $(1,\ldots ,1)$ (the relevant
vertex cannot be closed, and cannot have been previously
visited).  Therefore Odin cannot have a winning strategy,
so Eve wins.
\end{proof}

 As mentioned in the introduction, the method of this
section could likely be extended to obtain the conclusions of \cref{bcc} for
Trap on the Euclidean lattice $\Z^d$ for all $d$.  The main required step is
to show the analogue of \cref{p:frobose} for the following variant bootstrap
model defined on the odd vertices of $\Z^d$ (i.e.\ vertices whose coordinates
have odd sum):  if all but one of the $2d$ odd neighbors of any even vertex
are occupied, then the final odd neighbor becomes occupied at the next step.

We emphasize that the bootstrap argument of this section is not tight.  For
example, in the bottom-right picture of \cref{f:gameplay}, Eve wins although
the initial vertex does not become occupied in the bootstrap model.  To get
better bounds, we turn our attention to matchings.

\section{Matchings and Independent Sets}
\label{s:finite}

In this section we prove \cref{t:matchgame,t:indgame}.
Throughout this section $G=(V,E)$ will be a finite,
connected, simple, undirected graph. For $W\subseteq V$,
let $G\setminus W$ denote the subgraph of $G$ induced by
$V\setminus W$.

A \df{matching} $M$ of $G$ is a set of edges of $G$ no two
of which share a vertex. We say that $M$ \df{matches} a
vertex $v$, or that $v$ is \df{matched}, if $v$ is incident
to some edge of $M$.  In that case, the other incident
vertex to this edge is called the \df{partner} of $v$.  An
\df{independent set} is a set of vertices no two of which
are adjacent.  By a \df{maximum} matching or independent
set we mean one of maximum cardinality.

\begin{proof}[Proof of \cref{t:matchgame}]
In a game of Trap on $G$, if the first move is from $v$ to
$w$ then the remainder of the game is clearly equivalent to
Trap on $G\setminus \{v\}$ with initial vertex $w$.  We
will prove the claimed result by induction on the number of
vertices of $G$. In the base case $V=\{v\}$, the vertex $v$
is not in the maximum  matching, and indeed the first
player loses.

Suppose that $v$ lies in every maximum  matching of $G$,
and fix one such matching $M$.  Let $w$ be the partner of
$v$ in $M$, and let the first player move to $w$.  We claim
that $M':=M\setminus\{\{v,w\}\}$ is a maximum matching of
$G\setminus\{v\}$ (so that in particular $w$ does not lie
in every maximum  matching of this graph). Indeed, $M'$ is
clearly a maximum  matching of $G\setminus\{v,w\}$. Any
larger matching $M''$ of $G\setminus\{v\}$ must therefore
contain $w$.  But from such an $M''$, we could obtain a
maximum matching of $G$ with $v$ not matched by matching
$w$ to $v$ instead of its partner in $M''$.  This is a
contradiction, proving the claim. Now by the inductive
hypothesis, the next player loses on $G\setminus \{v\}$
starting from $w$.  Thus the first player wins on $G$
starting from $v$.

Now suppose that $v$ does not lie in every maximum matching
of $G$.  Suppose that the first player moves to any
neighbor $w$ of $v$.  We claim that $w$ lies in every
maximum  matching of $G\setminus\{v\}$.  A maximum matching
of $G$ that does not contain $v$ is also a maximum matching
of $G\setminus\{v\}$, so maximum matchings of $G$ and
$G\setminus\{v\}$ have the same size.  But if $M$ is a
maximum matching of $G\setminus\{v\}$ that does not contain
$w$, then $M\cup\{\{v,w\}\}$ is a larger matching of $G$, a
contradiction. This proves the claim.  Now by the inductive
hypothesis, the next player wins on $G\setminus \{v\}$
starting from $w$.  But $w$ was an arbitrary neighbor of
$v$, so the first player loses on $G$.
\end{proof}

\pagebreak
\begin{proof}[Proof of \cref{t:indgame}]
The proof is again by induction on the number of vertices. In the base case
$V=\{v\}$, the vertex $v$ lies in every maximum independent set, and indeed
the first player loses Vicious Trap.

Suppose that there exists a maximum independent set that does not contain
$v$.  Any such set must contain a neighbor of $v$, otherwise it could be
enlarged by adding $v$.  From among the maximum independent sets that do not
contain $v$, let $I$ be one that contains the fewest possible neighbours of
$v$.  Let the first player destroy the set $W$ of all neighbors of $v$ that
do not lie in $I$, and move to some other neighbor $w\in I$.  We claim that
every maximum independent set of $G\setminus (W\cup\{v\})$ contains $w$.
Indeed, $I$ is such a set, and any other independent set of the same size
would also be a maximum independent set of $G$ containing fewer neighbors of
$v$, contradicting the choice of $I$.  By the inductive hypothesis, the next
player loses.  Therefore, the first player wins.

Now suppose that $v$ lies in every maximum independent set of $G$.  Let $I$
be a maximum independent set of $G$.  Then $I':=I\setminus\{v\}$ is a maximum
independent set of $G\setminus \{v\}$. (Indeed, if $I''$ is a larger
independent set of $G\setminus \{v\}$ then either it contains no neighbor of
$v$, in which case it is an independent set of $G$ not containing $v$, or
else adding $v$ gives an independent set of $G$ that is larger than $I$).
Now, if the first player destroys a set $W$ of neighbors of $v$ and moves to
$w\not\in W$, then $I'$ is also a maximum independent set of $G\setminus
(W\cup\{v\})$. Since $I'$ does not contain $w$, the next player wins, by the
inductive hypothesis.
\end{proof}

\section{Upper bound: matching all odd vertices}
\label{s:ub}

In this section we prove \cref{upper}, which states that Eve typically wins
Trap on the diamond $D_n$ if $n$ is sufficiently large as a function of $p$.
This will be proved via \cref{t:matchgame}, by showing that there exist
appropriate matchings in $D_n$.

Throughout, we suppose that each odd vertex of the diamond $D_n$ is closed
with probability $p$, independently for different vertices, and all even
vertices are open.  We denote the associated probability measure $\P=\P_p$.
It will be convenient to consider arbitrary (partial) matchings of the
diamond $D_n$ itself, not just of its open subgraph.  Edges incident to
closed vertices in such a matching will be irrelevant in the eventual
application to the game.  Here is the key result of this section.

\begin{samepage}
\begin{proposition}
\label{p:ew} There exists $C>0$ and an event $\mathcal{E}
=\mathcal{E}_n$ such that if $p\rightarrow 0$ and
$n\rightarrow \infty$ in such a way that $n> (C\log
p^{-1})/p$ then we have $\P_{p}(\mathcal{E}_n)\rightarrow
1$, and such that on $\mathcal{E}_n$ we have the following.
\begin{enumerate}
\item[(i)] There exists a matching $M$ of $D_n$ that
    matches all open odd vertices.
\item[(ii)] For every protected even vertex $v$, there
    exists a matching $M_v$ in $D_n$ that matches all
    open odd vertices but leaves $v$ unmatched.
\end{enumerate}
\end{proposition}
\end{samepage}

Here is some motivation for \cref{p:ew}.  The diamond $D_n$ has $(2n)^2$ odd
vertices but only $(2n-1)^2$ even vertices.  Therefore, if we are to match
all open odd vertices as in (i), then at least $4n-1$ odd vertices must be
closed.  For this to happen with high probability we certainly require that
$p>1/n$. \cref{p:ew} states that it suffices to take $p$ larger than this
bound by a logarithmic factor, and that in that case we also get the stronger
conclusion (ii).  See \cref{p:enw} in the next section for a complementary
result in the other direction.

For the proofs it will be convenient to introduce an
alternative coordinate system.  We think of the diamond
$D_n$ as rotated 45 degrees clockwise, so that $(2n-1,0)$
is the bottom-right corner and $(-2n+1,0)$ is the top-left
corner. The following notation reflects this convention.
For $i=-2n+1, -2n+2, \ldots, 2n-1$, let $C_i=\{(x,y)\in
D_n:x+y=i\}$ be the \df{$i$th column} of the diamond.
Notice that if $i$ is odd then $C_i$ comprises $2n$ odd
vertices, while if $i$ is even then $C_i$ comprises $2n-1$
even vertices. We shall call the two cases {\bf odd
columns} and {\bf even columns} respectively. Similarly,
for $j=-2n+1,\ldots, 2n-1$, let $R_j=\{(x,y)\in
D_n:y-x=j\}$ be the \df{$j$th row} (which may again be odd
or even).  If $i$ and $j$ have the same parity, then $C_i$
and $R_j$ intersect in a unique vertex, which we write as
$$\langle i,j \rangle : = \Bigl(\frac{i-j}{2}, \frac{i+j}{2}\Bigr).$$
(If $i$ and $j$ have opposite parity, $C_i$ and $R_j$ do
not intersect). See \cref{f:rowcol}.
\begin{figure}
\begin{center}
\begin{tikzpicture}[scale=0.85,>=angle 90,scale=.8]
\def\n{3};
\begin{scope}[rotate=45]

\foreach \x in {0,...,3}
\path[fill=red!40] (.5-\x,1.5-\x)--++(1,0)--++(0,1)--++(-1,0)--cycle;
\foreach \x in {0,...,2}
\path[fill=blue!40] (-.5+\x,1.5-\x)--++(1,0)--++(0,1)--++(-1,0)--cycle;

\foreach \x in {0,...,\n} \draw
(\x+.5,-\n-.5+\x)--(\x+.5,\n+.5-\x)
(-\x-.5,-\n-.5+\x)--(-\x-.5,\n+.5-\x)
(-\n-.5+\x,\x+.5)--(\n+.5-\x,\x+.5)
(-\n-.5+\x,-\x-.5)--(\n+.5-\x,-\x-.5) ;

\draw[line width=2.5pt,green!40!black]
 (.5,.5)--++(1,0)--++(0,1)--++(1,0)
 --++(0,-1)--++(1,0)--++(0,-1)
 --++(-1,0)--++(0,-1)--++(-2,0)
 --cycle;

\coordinate (r)  at (0,2); \coordinate (c)   at (1,2);
\coordinate (o)  at (0,0); \coordinate (x)  at (0,-3);
\coordinate (q) at (2,-1);

\end{scope}

\path (r) -- ++(-1,0) coordinate (rl); \draw[
thick,<-,blue!50!black] (rl) -- ++(-2,0) node[left]
{$R_2$};

\path (q) -- ++(.8,0) coordinate (q1);
\draw[thick,<-,green!50!black] (q1) --++(1,0) node[right]
{$Q$};

\path (c) -- ++(0,.8) coordinate (cl); \draw[
thick,<-,red!50!black] (cl) to[out=90,in=180] ++(5,.5)
node[right] {$C_{-1}$};

\draw[thick,<-] (o) -- (-4.2,-.3) node[left] {$\langle
0,0\rangle$}; \draw[thick,<-] (x) -- ++(2,0)
node[right] {$\langle 3,-3\rangle$};

\end{tikzpicture}
\caption{The diamond $D_2$: coordinates, an even row, an
odd column,
 and the top-right quadrant.}
\label{f:rowcol}
\end{center}
\end{figure}

We divide the diamond into \df{quadrants}. The top-right
quadrant is given by
$$Q=Q^0=Q_n^0:=
\biggl(\bigcup_{i=0}^{2n-1}R_i\biggr) \cap
\biggl(\bigcup_{i=1}^{2n-1}C_i\biggr),$$
 that is, the vertices in top $2n$  rows and the
rightmost $2n-1$ columns.  For $k=1,2,3$ we define the
quadrant $Q^k=\theta^k(Q)$, where $\theta$ is the
anticlockwise rotation by 90 degrees about the origin.
Notice that the diamond can be written as the disjoint
union
$$D_n= \bigl\{\langle 0,0 \rangle\bigr\} \cup \bigcup_{k=0}^3 Q^k.$$

\sloppypar Given a matching $M$ in a graph, we say that a directed path
$(v_0,v_1,\ldots , v_{\ell})$ is \df{$M$-alternating} if it is self-avoiding,
and every other edge of the path starting with the first one belongs to $M$,
i.e., $\{v_i,v_{i+1}\}\in M$ for all even $i$. We will use the following
simple fact.  Recall that we allow matchings to include closed vertices.

\begin{lemma}
\label{l:alt}  Let $G$ be a finite bipartite graph in which
the odd vertices are declared open or closed. Let $M$ be a
matching that matches all open odd vertices. Let $v$ be an
even vertex that is matched in $M$. If there is an
$M$-alternating path starting from $v$ that contains a
closed odd vertex, then there exists a matching $M_v$ that
matches all open odd vertices but does not match $v$.
\end{lemma}

\begin{proof}
Let $(v_0,v_1,\ldots)$ be an $M$-alternating path starting
at $v=v_0$, and let $v_{2\ell+1}$ be the first closed odd
vertex on the path.  Construct $M_v$ from $M$ by removing
the edges $\{v_0,v_1\}, \{v_2,v_3\},\ldots ,\{v_{2\ell },
v_{2\ell +1}\}$ and adding the edges $\{v_1,v_2\},
\{v_3,v_4\},\ldots ,\{v_{2\ell -1}, v_{2\ell}\}$.
\end{proof}

The following is the main step of our proof of \cref{p:ew}.

\begin{samepage}
\begin{lemma}\label{l:ew2}
Let $\mathcal{F}=\mathcal{F}_n$ be the event that in each odd row $R_j\cap Q$
of the quadrant $Q$, there is a closed vertex $\langle i,j \rangle$ that is
not at either extreme end, i.e.\ $i\notin\{ 1, 2n-1\}$.  On $\mathcal{F}$,
there exists a matching $M$ of $Q$ with the following properties.
\begin{ilist}
\item
  All open odd vertices are matched.
\item  For each protected even vertex $v\in Q$, either:
\begin{alist}
  \item there exists an $M$-alternating path in $Q$ from
  $v$ to the top-left corner
    $\langle 1, 2n-1 \rangle$ of $Q$, or
  \item there exists an $M$-alternating path in $Q$,
  containing a closed vertex,  from
  $v$ to the top-right corner
    $\langle 2n-1,2n-1\rangle$ of $Q$.
\end{alist}
\item Each vertex $\langle 2n-1,j\rangle$ in the rightmost
    column $Q\cap C_{2n-1}$
    is matched in the down-left direction, to $\langle
    2n-2,j-1\rangle$.
\end{ilist}
\end{lemma}
\end{samepage}

\begin{figure*}
\begin{center}
\includegraphics[width=.7\textwidth]{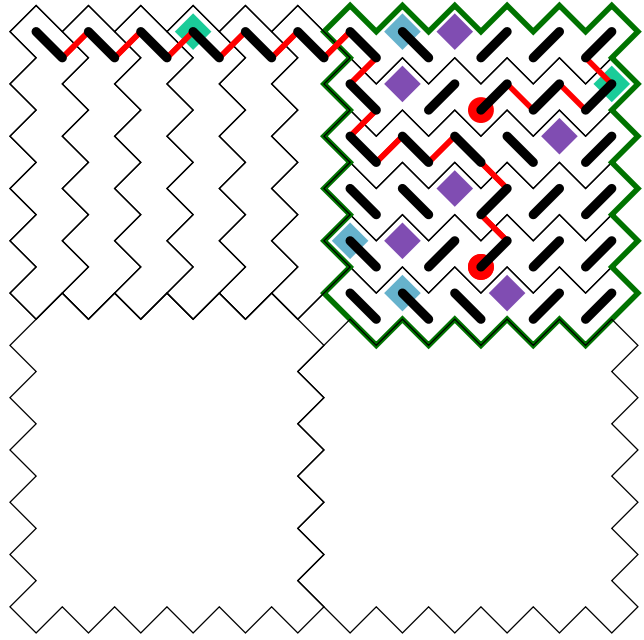}
\caption{Matching and alternating paths used in the proof of \cref{l:ew2}.
The matching $M$ in the top-right quadrant $Q$ of the diamond $D_6$ is shown
(together with the partition of the rows of $Q$ into pairs used in the
construction of $M$). Closed vertices are shown as filled squares, with
distinguished colors for those closed vertices comprising the set $H$, and
those in the rightmost column. Alternating paths from two protected vertices
(discs) are shown, ending at the top-right and top-left corners of the
quadrant. In the latter case, an extension into the top-left quadrant
 is also shown.  (Note that the alternating paths are used to
construct further matchings,
and are not directly related to trajectories of the token.)} \label{f:altpath}
\end{center}
\end{figure*}
\begin{proof}
The construction is illustrated in \cref{f:altpath}. On
$\mathcal{F}$, fix a set $H$ of closed vertices as follows.
For each $j\in \{1,3,\ldots, 2n-1\}$, let $m_j\in
\{3,5,\ldots ,2n-3\}$ be the largest number such that
$\langle m_j,j \rangle$ is closed, so that this is the
rightmost closed vertex in the $j$th row of $Q$ barring the
rightmost column. (This choice will be important in the
proof of property (ii)). Let $H=\bigl\{\langle m_j,j
\rangle : j\in \{1,\ldots , 2n-1\}\bigr\}$.  We construct
the matching $M$ in such a way that each vertex in an even
row is matched to some vertex in the odd row above it,
avoiding the vertices in $H$.  Specifically, for each $j\in
\{1, \ldots, 2n-1\}$, if $i>m_j$ then $(\langle i,j \rangle
, \langle i-1, j-1 \rangle) \in M$; if $i<m_j$ then
$(\langle i,j \rangle , \langle i+1, j-1 \rangle) \in M$.
Every odd vertex in $Q\setminus H$ is matched, so every
open odd vertex is matched (as well as every even vertex
and perhaps some closed odd vertices), so (i) holds.  Since
$H$ has no vertices in the rightmost column, (iii) holds.

We now proceed to check (ii). Fix a protected even vertex
$\langle i,j \rangle\in Q$. So there exists a closed odd
vertex $\langle i^*,j^* \rangle\in Q$ with $i^*> i$ and
$j^*> j$. We consider two cases, which will correspond
 to the two cases in the conclusion of part (ii).

\paragraph{\em Case (a)} {\em  There exists $\langle i^*,j^*
\rangle\in H$ with $i^*> i$ and $j^*> j$.}

Let $\langle i^*,j^*\rangle$ be the lowest element of $H$
that is above and right of $\langle i,j\rangle$, i.e.\ the
vertex satisfying the above condition for which $j^*$ is
smallest. For $j\leq k<j^*$, the even vertex $\langle
i,k\rangle$ is to the right of the element of $H$ in the
row immediately above it, so it is matched in the up-right
direction.  Therefore, there is an $M$-alternating path
$\pi_1$ from $\langle i,j \rangle$ to $\langle i,j^*-1
\rangle$ consisting of alternate up-right and up-left
steps.  (The path is empty if $j^*=j+1$).

Since $\langle i,j^*-1 \rangle$ is to the left of $\langle
i^*,j^*\rangle$, it is matched in the up-left direction.
Thus, there is an $M$-alternating path $\pi_2$ (again
possibly empty) from $\langle i,j^*-1 \rangle$ to $\langle
2,j^*-1 \rangle$ consisting of alternate up-left and
down-left steps.

Finally, since $H$ has no elements in the leftmost column
of $Q$, there is an $M$-alternating path $\pi_3$ from
$\langle 2,j^*-1 \rangle$ to $\langle 1,2n-1 \rangle$
consisting of alternate up-left and up-right steps.
Concatenating $\pi_1,\pi_2,\pi_3$ gives a path to the
top-left corner of the quadrant, as required.

\pagebreak
\paragraph{\em Case (b).} {\em There does not exist $\langle i^*,j^*
\rangle\in H$ with $i^*> i$ and $j^*> j$.}

Recall that in each odd row, the element of $H$ is the rightmost closed odd
vertex barring the rightmost column. Since $\langle i,j \rangle$ is
protected, it must therefore be protected by a closed vertex in the rightmost
column. Thus there exists $j^*>j$ with $\langle 2n-1,j^* \rangle$ closed.

Since $\langle i,j \rangle$ lies to the right of the
element of $H$ in the row immediately above it, there is an
$M$-alternating path from $\langle i,j \rangle$ to $\langle
2n-1,j \rangle$ using alternate up-right and down-right
steps, and thence to the top-right corner $\langle
2n-1,2n-1 \rangle$ using alternate up-right and up-left
steps.  This path passes through the closed vertex $\langle
2n-1,j^* \rangle$.
\end{proof}

\begin{proof}[Proof of Proposition \ref{p:ew}]
Let $\mathcal{F}=\mathcal{F}_n$ be the event in
\cref{l:ew2}.  Also let $\mathcal{G}=\mathcal{G}_n$ be the
event that the rightmost column $Q\cap C_{2n-1}$ of the
quadrant $Q$ contains a closed odd vertex.  For
$k=0,\ldots,3$, let $\mathcal{F}^k=\mathcal{F}^k_n$ and
$\mathcal{G}^k=\mathcal{G}^k_n$ be the images of
$\mathcal{F}$ and $\mathcal{G}$ under the rotation
$\theta^k$, i.e.\ the corresponding events in the rotated
quadrant $Q^k$.  We take
$$\mathcal{E}=\mathcal{E}_n:=\bigcap_{k=0}^3 \bigl(\mathcal{F}^k
\cap \mathcal{G}^k\bigr).$$

We first show that on $\mathcal{E}$, the properties (i) and (ii) in the
statement of the proposition hold.  The event $\mathcal{F}^k$ guarantees the
existence of a matching $M^k$ of the quadrant $Q^k$, the image under the
rotation $\theta^k$ of the matching in \cref{l:ew2}.  Let $M:=\bigcup_{k=0}^3
M^k$.  Since the origin is even, $M$ matches all odd vertices, as required
for (i).  Let $v\in D_n$ be an even vertex.  If $v$ is the origin then it is
unmatched in $M$, as required for (ii).  Otherwise, by symmetry we may assume
$v\in Q$.  By \cref{l:ew2}~(ii) there exists an $M$-alternating path from $v$
that either contains a closed vertex, or ends at the top-left corner of $Q$.
But in the latter case we can extend this path along the top edge of $Q^1$
using alternate down-left and up-left steps, by \cref{l:ew2} (iii), as in
\cref{f:altpath}. Since $\mathcal{G}^1$ holds, the resulting path then also
contains a closed vertex.  Now applying \cref{l:alt} gives a matching that
matches all open odd vertices but leaves $v$ unmatched, as required for (ii).

It remains to estimate $\P_{p}(\mathcal{E}_n)$. By a union
bound we get that $\P_p(\overline{\mathcal{F}_{n}}) \leq
n(1-p)^{n-2}$. Also, $\P_p(\overline{\mathcal{G}_{n}})=
(1-p)^{n}$. Therefore,
$$\P_p(\mathcal{F}_n\cap \mathcal{G}_{n}) \geq
1-(n+1)(1-p)^{n-2}\geq 1-2n e^{-pn/2},$$ provided $n\geq
4$.  For each fixed $p$, the expression on the right side
of the last inequality is increasing in $n$ for $n>2/p$.
Hence for fixed $C$, if $p$ is sufficiently small and  $n
> (C \log p^{-1})/p$, we have
$$\P_p(\mathcal{F}_n\cap\mathcal{G}_{n}) \geq
1-\frac{2C\log p^{-1}}{p} \exp \frac{-C\log p^{-1}}{2},$$
 which tends to $1$ as $p\to 0$ provided $C>2$.
By rotational symmetry, the probabilities
$\P_{p}(\mathcal{F}^k \cap \mathcal{G}^k)$ are equal for
$k=0,\ldots,3$, so another union bound shows that
$\P_p(\mathcal{E}_n)\to 1$.
\end{proof}

The proof of \cref{upper} is now straightforward.

\begin{proof}[Proof of \cref{upper}]
Choose $C$ as in Proposition \ref{p:ew}, and let $p\to 0$ and $n > (C\log
p^{-1})/p$. On $\mathcal{E}_n$, let $M$ be the matching of \cref{p:ew} (i).
Modify $M$ by removing from the matching each edge that is incident to a
closed vertex, to give a matching of the open subgraph of $D_n$. Since every
open odd vertex is matched, this is a maximum matching, and every maximum
matching matches all open odd vertices.  By \cref{t:matchgame}, we deduce
that Eve wins from every open odd site. Recall also that by convention Eve
also wins starting from a closed odd site.

Similarly, on $\mathcal{E}_n$, consider a protected even
vertex $v$, and the associated matching $M_v$ from
\cref{p:ew} (ii).  Removing edges incident to closed
vertices as before gives a maximum matching in which $v$ is
unmatched, so \cref{t:matchgame} gives that Eve wins from
$v$.
\end{proof}

\section{Lower bound: matching all even vertices}
\label{s:lb}

In this section we prove \cref{lower}. We adopt the same
coordinate system for the diamond $D_n$ as in the last
section, and we again consider the measure $\P_p$ under
which odd vertices of $D_n$ are closed with probability
$p$, and all even vertices are open. The following is the
main result in this section.

\begin{proposition}
\label{p:enw} Take any $c>0$ and let $n < c/(p \log
p^{-1})$ with $p\to 0$.  There exists an event
$\mathcal{O}=\mathcal{O}_n$ with $\P_p(\mathcal{O}_n)\to
1$, such that on $\mathcal{O}_n$, for every odd vertex
$v\in D_n$ there exists a matching $M_v$ of $D_n$ in which
all even vertices are matched, but neither $v$ nor any
closed odd vertex is matched.
\end{proposition}

The proof of \cref{p:enw} will use the next lemma.  In the
application, the set $H$ will consist of the closed
vertices together with one additional arbitrary odd vertex
$v$.  The lemma will in fact be applied twice: both to rows
and to columns.  By an interval of $s$ consecutive rows we
mean a set of the form $\bigcup_{j=1}^s R_{j+a}$, and by an
interval of $s$ consecutive odd rows we mean a set of the
form $\bigcup_{j=1}^{s} R_{2j+2a+1}$, where $a\in\Z$.

\begin{samepage}
\begin{lemma} \label{l:enw} Let $n$ and $s$ be positive integers.  Let
$H$ be a set of odd vertices of the diamond $D_n$.  Suppose that every
interval of $s$ consecutive odd rows
contains at most $s$ vertices of $H$, and that for each
$\langle i,j \rangle \in H$, there is no other $\langle
i,j' \rangle \in H$ in the same column with $|j-j'|< 2s$.
Then there exists a matching of $D_n$ that matches all even
vertices but leaves $H$ unmatched.
\end{lemma}
\end{samepage}

Note that in the above lemma we do not require that $s\leq
2n$. If $s>2n$, then the first condition in the lemma is
vacuously satisfied, and the second condition states that
no odd column contains two or more elements of $H$.  In the
main case of interest, we will in fact choose $s$ to be
approximately $\log p^{-1}/\log\log p^{-1}$. \cref{l:enw}
will in turn be proved by partitioning the rows into
suitable intervals, and using the following technical
lemma.
\begin{samepage}
\begin{lemma}\label{l:enwind}
Fix $n$, and consider an interval of rows of the diamond
$$W:=\bigcup_{j=a}^b R_j \subseteq D_n,$$
where $a\leq b$ and $a$ is odd.  Let $H$ be a subset of the
odd vertices of $W$, such that each column contains at most
one vertex of $H$.  Suppose that either:
\begin{ilist}
\item $b$ is odd, or
\item $b$ is even, and for each integer $t$, the top $t$
    odd rows of $W$ contain at most $t$ vertices of $H$.
\end{ilist}
Then there exists a matching of $W$ that matches all even
vertices but leaves $H$ unmatched.
\end{lemma}
\end{samepage}

\begin{figure*}
\begin{center}
\includegraphics[width=0.55\textwidth]{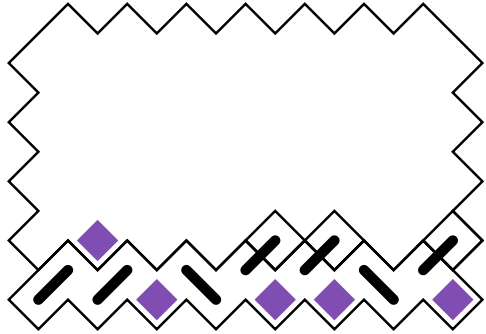}
\caption{Inductive step in the proof of \cref{l:enwind} (in case (ii)).
Vertices of $H$ are shown as filled squares.  The inductive hypothesis is
applied after removing the bottom two rows, and modifying $H$ by adding the
three outlined vertices.} \label{f:lb}
\end{center}
\end{figure*}
\begin{proof}  See \cref{f:lb}.
The proof is by induction on $s:=\lfloor (b-a)/2\rfloor$,
the number of even rows of $W$.  The inductive step will
involve removing the bottom two rows of $W$ and modifying
$H$. When $s=0$, either $W$ consists of one odd row (case
(i)), or $W$ is empty (case (ii)); in either case we take
the empty matching.

Now suppose $s\geq 1$. Let $z=\langle \ell, a\rangle$ be
the leftmost vertex of $H$ in the bottom row $R_a$.  (If
none exists, take $\ell=\infty$.)  Consider any even site
$u=\langle i,a+1\rangle$ in the next row $R_{a+1}$.  If $u$
is to the left of $z$ (i.e.\ $i<\ell$), match it down-left
to $\langle i-1,a\rangle$.  If $u$ is to the right of $z$
(i.e.\ $i>\ell$), match it down-right to $\langle
i+1,a\rangle$, unless the latter vertex is in $H$.  In that
case, match $u$ up-right to $\langle i+1,a+2\rangle$
instead. Since $H$ has at most one vertex in each column,
this last vertex is not in $H$.

Let $W':=\bigcup_{j=a+2}^b R_j$, and let
$$H':=\bigl(H\cap W'\bigr)\cup \Bigl\{\langle i+1,a+2\rangle :
 i>\ell \text{ and } \langle
i+1,a\rangle \in H\Bigr\}.$$
 (Thus $H\cap W'$ is augmented by the set of odd vertices of
$R_{a+2}$ that were already matched to vertices in
$R_{a+1}$; this corresponds to shifting the vertices of $H$
in the bottom row up by $2$, except for the leftmost one).
We will apply the inductive hypothesis to $W'$ and $H'$,
and combine the resulting matching with the matching
constructed above.

To complete the argument, we must check that $W'$ and $H'$
indeed satisfy the conditions of the lemma.  Since the
construction of $H'$ involves shifting vertices of $H$
vertically, $H'$ still has no two vertices in the same
column. If (i) holds for $W$ then (i) obviously holds for
$W'$ as well. If (ii) holds for $H$, we need only check
condition (ii) for $H'$ in the case $t=s-1$, i.e.\ that
$|H'|\leq s-1$. If $\ell=\infty$ then $H'=H\cap W'$ and so
this follows immediately from the condition on $H$. If
$\ell<\infty$, then we have $|H|\leq s$, but the shifted
vertex $\langle \ell,a+2\rangle$ corresponding to $z$ is
{\em not} included in $H'$, so $|H'|=|H|-1\leq s-1$.
\end{proof}

\begin{proof}[Proof of \cref{l:enw}]
We partition the rows of the diamond, starting from the
bottom, into minimal intervals containing no more vertices
of $H$ than odd rows. More precisely, let $\ell_0=-2n+1$,
and iteratively define $\ell_{k+1}$ to be the smallest odd
integer in $(\ell_k,2n-1]$ for which
$$\Bigl| H \cap {\textstyle \bigcup_{j=\ell_k}^{\ell_{k+1}-1} R_j} \Bigr| \leq
\frac{\ell_{k+1}-\ell_k}{2};$$
 if there is no such number then we instead take
 $\ell_{k+1}=2n$, write $K:=k$, and stop the iteration.

Since any interval of $s$ rows contains at most $s$
vertices of $H$, each of these intervals must contain at
most $s$ odd rows (i.e.\ $\ell_{k+1}-\ell_k\leq 2s$), and
therefore each of the corresponding regions
$W_k:=\bigcup_{j=\ell_k}^{\ell_{k+1}-1} R_j$ contains no
two vertices of $H$ in the same column. Furthermore, in
each of these regions except perhaps the last, $W_K$, for
every $t$, the top $t$ odd rows contain at most $t$
vertices of $H$ (otherwise $\ell_{k+1}$ should have been
smaller). On the other hand, the last region $W_K$ has an
odd row at the top. Therefore we can apply \cref{l:enwind}
to each of the regions $W_0,\ldots,W_K$ and take the union
of the resulting matchings.  (If $s\geq 2n$, it is possible
that $K=0$ and thus $W_K=D_n$.)
\end{proof}

\begin{proof}[Proof of Proposition \ref{p:enw}]
Fix any $c>0$. Let
\begin{equation}
\label{e:choosek}
s=\biggl\lceil\frac{4\log p^{-1}}{\log\bigl[(\log p^{-1})/4c\bigr]}\biggr\rceil.
\end{equation}
We define the following events. Let
$\mathcal{R}=\mathcal{R}_{n,s}$ be the event that in any
$s$ consecutive odd rows of $D_n$, at most $s-1$ vertices
are closed. Let $\mathcal{T}=\mathcal{T}_{n,s}$ be the
event that there are no two distinct closed vertices
$\langle i, j \rangle$ and $\langle i,j' \rangle$ in the
same column with $|j-j'| <2s$.
Finally, for an odd vertex $v=\langle i,j \rangle$ we
define $\mathcal{X}_{v}=\mathcal{X}_{v,n,s}$ to be the
event that there is no closed vertex $\langle i,j' \rangle$
with $|j-j'|< 2s$ and $j\neq j'$.

Note that on $\mathcal{R}\cap
\mathcal{T}\cap\mathcal{X}_v$, the hypothesis of
\cref{l:enw} is satisfied when $H$ is taken to be the set
of closed vertices together with $v$.  (The addition of $v$
is the reason for the using $s-1$ in definition of
$\mathcal{R}$). Thus there exists a matching $M_v$
satisfying the required conclusion for $v$. However, we
cannot directly obtain the same conclusion simultaneously
for all $v$: indeed, the event $\mathcal{X}_v$ does not
hold for vertices $v$ within distance $2s$ of a closed
vertex in the same column.

To address this issue we also consider rotated versions of
the same events.  Let $\mathcal{R}',\mathcal{T}'$ be the
images of $\mathcal{R},\mathcal{T}$ under the $90$ degree
anti-clockwise rotation $\theta$ about the origin. Also for
an odd vertex $v=\langle i,j \rangle$, let $\mathcal{X}'_v$
be the event that there is no closed vertex $\langle i',j
\rangle$ with $|i-i'|<2s$ and $i\neq i'$. By rotational
symmetry, $v$ satisfies the required conclusion on
 $\mathcal{R}'\cap \mathcal{T}'
\cap \mathcal{X}'_v$ as well. It follows that on
$$\mathcal{O}=\mathcal{O}_n:=
\mathcal{R}\cap \mathcal{R}'\cap \mathcal{T}\cap
\mathcal{T}' \cap \bigcap_{v\text{ odd}} (
\mathcal{X}_v\cup \mathcal{X}'_v \bigr),$$ the required
conclusion holds for all odd vertices $v\in D_n$.

Now we estimate $\P_p(\mathcal{O}_n)$. For the following
calculations we always take $p$ to be sufficiently small
such that $p^{-1}>c^{-1}\log p^{-1}/4>5$ and $s< c/(p\log
p^{-1})$. Notice that this ensures that for all $n$ with
$1\leq n< c/(p \log p^{-1})$ we have $np< 1/20$, and that
$s$ as given by (\ref{e:choosek}) is at least $2$.

If $s>2n$ then $\P_p(\mathcal{R}_{n,s})=1$. Otherwise, a
union bound over all intervals of $s$ consecutive odd rows
gives
$$
\P_p\bigl(\mathcal{R}_{n,s}\bigr)\geq  1-2n\,\P\bigl(Z\geq s-1\bigr)
\geq  1-2n\,\P\bigl(Z\geq s/2\bigr),
$$
where $Z$ is a binomial random variable with parameters
$(2ns,p)$. Using $np<1/20$, a Chernoff bound then yields
\begin{equation}
\label{e:enw1}
\P_p(\mathcal{R}_{n,s}) \geq  1-2n(4np)^{s/4}.
\end{equation}
 The expression on the right is
decreasing in $n$ for fixed $p$ and $s> 0$.  Hence for $n <
c/(p \log p^{-1})$, using \eqref{e:choosek} we have
\begin{equation}
\label{e:enw4}
\P_p(\mathcal{R}_{n,s})
\geq 1-\frac{2c}{p \log p^{-1}}
\biggl(\frac{4c}{\log p^{-1}}\biggr)^{s/4}
\geq 1-\frac{2c}{\log p^{-1}}.
\end{equation}

Now note that in any given odd column, the probability that
there are two closed odd vertices with their vertical
coordinates differing by less than $2s$ is at most
$2nsp^2$. Taking a union bound over all odd columns, we
obtain for $n < c/(p \log p^{-1})$,
\begin{equation}
\label{e:enw2}
\P_p(\mathcal{T}_{n,s})\geq 1-4n^2sp^2\geq  1-\frac{4c^2s}{\log^2p^{-1}}.
\end{equation}

Also, for each odd vertex $v$ we have
$\P_p(\mathcal{X}_v\cup \mathcal{X}'_v)\geq 1-(2sp)^2$, and
a union bound then gives that for $n < c/(p \log p^{-1})$,
\begin{equation}
\label{e:enw3}
\P_p\biggl(\,\bigcap_{v~\text{odd}} (\mathcal{X}_v\cup \mathcal{X}'_v)\biggr)
\geq  1-4n^2(2sp)^2\geq 1-\frac{16c^2s^2}{\log^2p^{-1}}.
\end{equation}

Using the definition of $s$ in \eqref{e:choosek},  we see
that the right sides of \eqref{e:enw4}, \eqref{e:enw2},
\eqref{e:enw3} each converge to $1$ as $p\rightarrow 0$. By
rotational symmetry, we have
$\P_p(\mathcal{R}')=\P_p(\mathcal{R})$ and
$\P_p(\mathcal{T}')=\P_p(\mathcal{T})$. A final union bound
now shows that $\P_p(\mathcal{O}_n)\rightarrow 1$ as
$p\rightarrow 0$ with $n < c/(p \log p^{-1})$.
\end{proof}

\begin{proof}[Proof of \cref{lower}]
We apply \cref{p:enw}. On the event $\mathcal{O}_n$, for each odd $v$, the
matching $M_v$ is a maximum matching of the open subgraph of $D_n$, in which
all even vertices are matched but $v$ is unmatched.  By \cref{t:matchgame},
Odin wins from every initial vertex, whether even or odd.
\end{proof}

\section{Further results}
\label{s:game} In this section we justify the claims about
protected vertices in the introduction, and briefly address
some consequences of \cref{lower,upper} for Trap on finite
regions other than the diamond, for the density of closed
even vertices, and for the length of the game on $\Z^2$.

\subsection*{Protected vertices}
We check the remarks following \cref{upper} regarding the
set $S$. Write $L=C'\log p^{-1}/p$, as in the definition of
$S$ in \eqref{hyperbola}.  We adopt the rotated coordinate
system of \S~\ref{s:ub}.  If $\langle i,j\rangle\in D_n$
satisfies $i,j\geq 0$ and $(2n-i)(2n-j)>L$ then the
probability that there is a closed vertex of $D_n$ above
and right of $\langle i,j\rangle$ is at least
$$1-(1-p)^L\geq 1-e^{-pL}=1-p^{C'}.$$
To ensure that all vertices in $S$ are protected, we need
only check this condition (and similar ones involving other
quadrants) for $O(L)$ vertices near the boundary of $S$. (A
key point to note is that provided $2n>L+2$, one of the
vertices to be checked will be in the top row, of the form
 $\langle
2n-L-\epsilon,2n-1\rangle$ for some $\epsilon\in(0,2]$.)
 Therefore for $C'>1$,
a union bound shows that all vertices in $S$ are protected
with high probability as $p\to 0$, uniformly in $n$.

On the other hand, we have as $p\to 0$, uniformly in $n$,
$$|D_n\setminus S|=O\biggl(\int_1^L
\frac{L}{x}\,dx\biggr)=O(L\log L)=O\bigl(p^{-1} \log^2 p^{-1}
\bigr),$$ as claimed.

\subsection*{Finite regions}

Fix a finite connected region $A\subseteq \Z^2$ whose
internal boundary consists entirely of odd vertices.
Consider a game of Trap on the open subgraph of $A$, where
as usual each odd vertex is closed with probability $p$.
(Equivalently, we can consider Trap on $\Z^2$ started from
a vertex in $A$, but declaring a win for Odin if Eve ever
leaves $A$).  We can infer outcomes of this this game in
certain situations by comparing with the games played on
smaller and larger diamonds and using \cref{lower,upper}.
Suppose that $D_n'\subseteq A \subseteq D_N'$, where $D_n'$
and $D_N'$ are translates of diamonds $D_n$ and $D_N$.
Starting from any vertex in $A$, if Odin wins the game
played on $D_N$ then he wins the game on $A$ as well.  On
the other hand, if Eve wins in $D_n$ starting from a vertex
in $D_n$ then she wins on $A$.  However, the latter
argument does not apply to initial vertices in $A\setminus
D_n$.

As a concrete example, consider a square modified so that all internal
boundary vertices are odd:
$$B^\#(n):=\Bigl([1,n]^2\cap \Z^2\Bigr) \cup
\Bigl([0,n+1]^2 \cap \Z_{\rm o}^2\Bigr),$$
 where $\Z_{\rm o}^2$ is the set of odd vertices of $\Z^2$.  In the
following, $f(p)\gg g(p)$ means $f(p)/g(p)\to\infty$.
\begin{corollary}
Consider Trap on the open subset of $B^\#(n)$, where each
odd vertex is closed with probability $p$ and all even
vertices are open. If $p\to 0$ with $n< c/(p\log p^{-1})$
and any $c$, then Odin wins from every vertex with high
probability.  If $p\to 0$ with $n\gg \log p^{-1}/p$ then
with high probability, Eve wins from a proportion $1-o(1)$
of initial vertices.
\end{corollary}

\begin{proof}
Since $B^\#(n)$ is contained in a translate of $D_{n+1}$,
the first claim follows from \cref{lower}.  For the second
claim, let $m=m(p)=2C \log p^{-1}/p$ where $C$ is the
constant of \cref{upper}.  We can `approximately tile'
$B^\#(n)$ with disjoint translates of $D_m$ in such a way
that all but a proportion $o(1)$ of vertices are covered.
By the law of large numbers, with high probability, the
conclusion of \cref{upper} holds for a proportion $1-o(1)$
of these translates. Using the remarks following
\cref{upper} about protected vertices, the required
conclusion follows.
\end{proof}

\subsection*{Density of closed even vertices}

Recall that in dimension $2$, the body-centered lattice $\B^2$ and the square
lattice $\Z^2$ are isomorphic.  As mentioned in the introduction, it is
possible to give an alternative proof of \cref{bcc} for $d=2$ using
Proposition \ref{p:ew}, and this gives a much tighter bound on $q(p,2)$ in
part (ii) compared with the argument in \S~\ref{s:bcc}.

\begin{corollary}
Let odd and even vertices of $\Z^2$ be closed with respective probabilities
$p$ and $q$.  There exists $c_0>0$ such that if $p>0$ and $q< c_0 p^2
\log^{-2} p^{-1}$, then almost surely, from every initial vertex, Trap is not
a draw.
\end{corollary}

\begin{proof}
For an even vertex $u\in \Z^2$, call a diamond-shaped
region $D_n(u):=u+D_n$ \df{good} if it contains no closed
even vertices, and it has a matching that matches all open
odd vertices. Starting from any odd vertex in a good
diamond $D_n(u)$, Eve can win within $D_n(u)$. We now
proceed as in the proof of \cref{bcc}.  (A diamond $D_n(u)$
induces a graph in $\Z^2$ isomorphic to that induced by a
box $\widetilde{B}(u',n')$ in $\B^2$ as in \S~\ref{s:bcc}).
 It follows from \cref{p:ew} that
if $q < c_0 p^2\log^{-2} p^{-1}$ for a suitable constant
$c_0>0$, then choosing $n$ appropriately, the probability
that $D_n(u)$ is good exceeds one minus the critical
probability of the two-dimensional star-lattice. Almost
surely we can then surround the initial vertex by good
diamonds, and the proof goes through as before.
\end{proof}

\subsection*{Length of the game}

We can use \cref{upper,lower} to obtain bounds on the time
for Trap to terminate.  For simplicity we do this only for
$q=0$, although analogous results are also available in the
regime of the last result.

If Eve can win, then it is natural for her to try to win as
quickly as possible, while Odin tries to prolong the game.
Fix an initial vertex $v$, and let $\mathbf E$ and $\mathbf
O$ denote the sets of all possible strategies for Eve and
Odin respectively starting from $v$.  For $\mathbf{e}\in
\mathbf E$ and $\mathbf{o}\in\mathbf{O}$, let
$T(\mathbf{e},\mathbf{o})$ be the number of turns until the
game terminates when Eve plays with strategy $\mathbf e$
and Odin plays with strategy $\mathbf o$, provided Eve
wins; if Eve does not win with this pair of strategies, let
$T(\mathbf{e},\mathbf{o})=\infty$.  Define
$$T=T_v:=\adjustlimits\inf_{\mathbf{e}\in
\mathbf E}\sup_{\mathbf{o}\in \mathbf O} T(\mathbf{e},\mathbf{o}),$$
i.e.\
the minimum time in which Eve can guarantee to win.
\begin{corollary}
\label{cor:time} Let odd and even vertices of $\Z^2$ be
closed with respective probabilities $p>0$ and $q=0$. There
exist constants $c_1,C_1>0$ such that, for any initial
vertex $v$, with high probability as $p\to 0$, the time
$T=T_v$ for Eve to win Trap satisfies
$$\frac{c_1}{p\log p^{-1}} \leq T \leq \frac{C_1 \log^2 p^{-1}}{p^2}.$$
\end{corollary}

\begin{proof}
By translation-invariance we can assume without loss of generality that $v$
is $(0,0)$ or $(1,0)$.  Let $n=\lfloor c/(p\log p^{-1})\rfloor$ and $N=\lceil
(C \log p^{-1})/p\rceil$, where $c,C$ are the constants of
\cref{upper,lower}.  With high probability, Eve can guarantee a win without
having to leave $D_N$, and hence she can win in at most $|D_N|$ moves. On the
other hand, with high probability, Eve cannot win within $D_n$, and it takes
at least $2n-2$ steps for the token to leave $D_n$ starting from $v$.
\end{proof}

\section*{Open problems}

\begin{ilist}
  \item
  On $\Z^d$ with each vertex closed independently with
  probability $p$, does there exist $p>0$ for which Trap
  starting from the origin is a draw with positive
  probability?  (It is plausible that the answer is no for
  $d=2$ and yes for $d\geq 3$).
  \item
  In the situation of (i) above, is the probability of a draw
  monotone in $p$?
  \item
  On a diamond $D_n=\{u\in\Z^2: \|u\|_1<2n\}$, with odd vertices closed with
  probability $p$ and all even vertices open, what happens within
  the window $c/(\log p^{-1})<np<C \log p^{-1}$?  Is there
  a regime in which both players have substantial regions
  of winning initial vertices?
  \item
  On $\Z^2$ with odd vertices closed with
  probability $p$ and all even vertices open, what more can be said
  about the minimum time $T$ in which Eve can guarantee a
  win?  (\cref{cor:time} states that with high probability it is between
  $p^{-2}$ and $p^{-1}$ up to logarithmic factors).  Does
  $T$ converge when suitably scaled as $p\to 0$?
  \item
  How do our results change on non-bipartite graphs such as
  the triangular lattice?  How do the outcomes
  of Trap and Vicious Trap differ from each other?
  \item
  How do our results change for the mis\`ere variant of Trap in which a
  player who cannot move wins, or for the variant game in which one
  specified player wins if either player cannot move?  (See \cite{game-tree}
  for analysis of these games on random trees).
\end{ilist}

\section*{Acknowledgements}

R.B.\ was supported by a U.C.\ Berkeley Graduate Fellowship. J.B.M.\ was
supported by EPSRC Fellowship EP/E060730/1.  J.W.\ was supported by the
Swedish Research Council, the Knut and Alice Wallenberg Foundation, and the
G{\"o}ran Gustafsson Foundation.  Much of this work was completed at the
Theory Group of Microsoft Research in Redmond during visits by J.B.M.\ and
J.W.\ and an internship by R.B.  We thank Microsoft for its hospitality and
support.

\end{document}